\documentclass{amsart}
\usepackage[marginratio=1:1]{geometry}

\usepackage[hidelinks]{hyperref}
\usepackage{calc}
\newsavebox\CBoxx
\newcommand\hcancel[2][0.5pt]{%
  \ifmmode\sbox\CBox{$#2$}\else\sbox\CBox{#2}\fi%
  \makebox[0pt][l]{\usebox\CBox}%
  \rule[0.5\ht\CBox-#1/2]{\wd\CBox}{#1}}
\usepackage{blindtext}
\usepackage{color}
\usepackage{hyperref,url}
\usepackage{float}
\usepackage{hyperref}
\usepackage{enumerate}
\usepackage{enumitem}   
\usepackage{bbm}
\usepackage{cancel}
\usepackage{mathrsfs}
\usepackage{colonequals}
\usepackage{pdfpages}
\newcommand\smallO{
  \mathchoice
    {{\scriptstyle\mathcal{O}}}
    {{\scriptstyle\mathcal{O}}}
    {{\scriptscriptstyle\mathcal{O}}}
    {\scalebox{.7}{$\scriptscriptstyle\mathcal{O}$}}
  }

\usepackage{amsmath,amsfonts,amsthm,bm}
\usepackage{mathrsfs}  
\usepackage{xcolor}
\usepackage{amsfonts}
\usepackage{amssymb}
\usepackage{amsmath}
\usepackage{mathtools}
\usepackage{mathrsfs}  

\usepackage[normalem]{ulem}

\numberwithin{equation}{section}
\usepackage[toc,page]{appendix}
\usepackage{tikz-cd}
\theoremstyle{definition}

\newtheorem{definicao}{Definition}[section]
\newtheorem{remark}[definicao]{Remark}

\theoremstyle{plain}

\newtheorem{teorema}[definicao]{Theorem}

\newtheorem{lema}[definicao]{Lemma}
\newtheorem{proposicao}[definicao]{Proposition}
\newtheorem{corolario}[definicao]{Corollary}

\newtheorem {mtheorem} {\bf Theorem}

\newtheorem {mcorollary}[mtheorem]  {\bf Corollary}

\newtheorem{manualtheoreminner}{\bf Claim}

\newtheorem{manualtheoreminnerb}{\bf Step}
\newenvironment{step}[1]{%
  \renewcommand\themanualtheoreminnerb{#1}%
  \manualtheoreminnerb
}{\endmanualtheoreminner}

\newenvironment{theoremproofA}[1]{\par\noindent{\emph{Proof of Theorem A.}} \space#1}{\leavevmode\unskip\penalty9999 \hbox{}\nobreak\hfill\quad\hbox{$\qed$}}
\newenvironment{theoremproofB}[1]{\par\noindent{\emph{Proof of Corollaries B and C.}}  }{\leavevmode\unskip\penalty9999 \newline\hbox{}\nobreak\hfill \quad\hbox{$\qed$}}

\usepackage{scalerel}

\definecolor{roxo}{rgb}{0.44, 0.16, 0.39}
\definecolor{ao(english)}{rgb}{0.0, 0.5, 0.0}
\definecolor{dmagenta}{RGB}{139, 0, 139}
\definecolor{dgreen}{RGB}{0,90,0}
\definecolor{navy}{RGB}{0,0,128}

\DeclareSymbolFont{bbold}{U}{bbold}{m}{n}
\DeclareSymbolFontAlphabet{\mathbbold}{bbold}

\usepackage{stackengine}

\def\R{\mathbb R}
\def\d{\mathrm d}
\def \d {\mathrm{d}}

\def \tri {\mathbbold{\Delta}}

\def \dist{\mathrm{dist}}

\DeclareMathOperator*{\argmin}{arg\,min}

\definecolor{iblue}{RGB}{0, 35, 194}
\title[Random Young towers for predominantly expanding multimodal circle maps]{Random Young towers and quenched decay of correlations for predominantly expanding multimodal circle maps}
\author[M.M. Castro, G. Tenaglia.]
{Matheus M. Castro$^{*}$, Giuseppe Tenaglia$^{\dagger}$.
}

\address{${}^*$Anita B. Lawrence Centre (H13), University Mall, UNSW, Kensington NSW 2033.}

\address{${}^\dagger$Department of Mathematics, Imperial College London, London, SW7 2AZ, United Kingdom.}

\email{m.manzatto\_de\_castro@unsw.edu.au, giuseppe.tenaglia20@imperial.ac.uk.}

\allowdisplaybreaks
\begin{document}



\maketitle
\vspace{-1.3cm}
\begin{align*}
   {}^* \textit{\small University of New South Wales},\ {}^\dagger \textit{\small Imperial College London} 
\end{align*}
\begin{abstract}
In this paper, we study the random dynamical system $f_\omega^n$  generated by a family of maps $\{f_{\omega_0}: \mathbb{S}^1 \to \mathbb{S}^1\}_{\omega_0 \in [-\varepsilon,\varepsilon]},$ $f_{\omega_0}(x) = \alpha \xi (x+\omega_0) +a\ (\mathrm{mod }\ 1),$ where $\xi: \mathbb S^1 \to \mathbb R$ is a non-degenerated map, $a\in [0,1)$, and $\alpha,\varepsilon>0$. Fixing a constant $c\in (0,1)$, we show that for $\alpha$ sufficiently large and for $\varepsilon > \alpha^{c-1},$ the random dynamical system $f_\omega^n$ presents a random Young tower structure and quenched decay of correlations.

\end{abstract}

\section{Introduction}

The study of chaotic systems is often focused on understanding the statistical properties that result from their intricate dynamics. Within the class of chaotic dynamical systems, mixing dynamical systems hold particular importance. Specifically, a map $f:M\to M$, acting on a compact manifold $M$ with an invariant probability measure $\mu$ on $M$ (i.e., $\mu(\cdot) = \mu(f^{-1}(\cdot))$), is referred to as mixing if, for every bounded and measurable observables  $\psi, \varphi: M \to \mathbb{R},$ \begin{align}
   \left|\int_M \varphi \circ f^n \ \psi \, \d \mu - \int_M \varphi \,\d \mu \int_M \psi \, \d \mu\right| \xrightarrow[]{n\to\infty} 0. \label{correlation}
\end{align}
The equation above suggests that if $f$ is mixing, then the random variables $\varphi \circ f^n$ and $\psi$ tend towards independence.

When studying mixing systems, a natural question is how fast the convergence of \eqref{correlation} occurs. Two commonly employed techniques to address this question are the spectral gap method and the construction of Young towers for the system $f$. The first technique involves identifying a suitable Banach space $F\subset L^1(M,\mu)$ such that the Perron-Frobenius operator induced by $f$ exhibits the spectral gap property (see  \cite{CarTransf, BTransf,Btranf2,CarlangeloTransf,Ruelle1989} and the references therein). The second method was introduced by L.-S. Young in her series of seminal works {\text{
blue}\cite{LSYoung1,LSYoung2,LSYoung3,LSYoung,LSY}}. This approach involves identifying a subset $\Delta \subset M$ and defining a return function $R: \Delta \to \mathbb{N}$ such that the $\mu$-measurable map $f^R: \Delta \to \Delta$, given by $f^R(x) = f^{R(x)}(x),$ exhibits strong hyperbolic properties. The hyperbolic behaviour of $f^R$, combined with the decay rate of the sequence $\{\mu(R > n)\}_{n \in \mathbb{N}}$, provides the means to establish an upper bound for the speed of convergence of \eqref{correlation} (see \cite[Chapters 3 and 6]{Alves} for further details).

The question of how fast the convergence of \eqref{correlation} occurs can be extended to the random dynamical systems context. Given a random dynamical system $f_\omega^n: M\to M$ over a $\mathbb P$-invariant ergodic dynamical system $\theta:\Omega \to \Omega$ one may consider how fast
\begin{align}
   \left|\int_M \varphi \circ f_\omega^n \ \psi \, \d \mu_\omega - \int_M \varphi \,\d \mu_{\theta^n\omega} \int_M \psi \, \d \mu_{\omega}\right| \xrightarrow[]{n\to\infty} 0, \label{correlation1}
\end{align}
 where $\mu(\d x, \d \mu) = \mu_\omega(\d x) \mathbb P(\d \omega)$ is an invariant measure for the skew product $\Theta:\Omega\times M\to \Omega\times M$, $\Theta(\omega,x) = (\theta\omega,f_\omega(x)).$ In the literature \eqref{correlation1} is generally called \emph{quenched decay of correlation}.

The two aforementioned techniques have their random dynamical system counterparts. The first approach is a functional analytical method that depends either on the analysis of non-autonomous Perron-Frobenius operators (see \cite{Atnip, Buzzi,Froyland, VarTher}) or on the description of the stationary measures of the two-point process $f^{(2)}_\omega(x,y) =(f_\omega(x), f_\omega(y))$ 
 on $M\times M$ (see \cite{LToral, BS,Dolg, BMixing}). The second technique derives from the theory of random Young towers, which has been actively studied in recent years (see \cite{Alves1,Alves2,Baladi, Du, Su, LSY} and references therein).

The construction of random Young towers has proven challenging in several contexts. For transitive non-uniformly expanding maps, the construction of a Young tower has been established through two key observations. Firstly, given an arbitrarily small open set $U$, there exists a subset $V \subset U$ and an integer $n = n(U)$ such that $f^n(V)$ is a ball of fixed radius $\delta > 0$ and $f^n|_{V}$ is uniformly expanding. Secondly, there exists a suitable reference set $\Delta_0 \subset M$ and a natural number $L \in \mathbb{N}$ such that $f^{n+i}(V)$ contains $\Delta_0$ for some $i \in \{0, 1, \ldots, L\}$. Additionally, $f^{n+i}|_{V}$ is uniformly hyperbolic. These two features combined allow the construction of a Young tower \cite{Alves,Alves1}.
  
The above scheme has been successfully adapted to the random context in \cite{Araujo2003,Vil2018, Su}  to build a random Young tower for small perturbations of transitive non-uniformly expanding maps. Additionally, in \cite{Baladi} and \cite{Wael2019}, this approach has been combined with prior knowledge of mixing and hyperbolic properties of the systems under consideration. This allowed for the construction of random Young towers for small perturbations of the logistic family (within the range of Benedicks and Carleson parameters) and random iterations of LSV maps, respectively.

In this paper, we construct random Young towers and show their quenched decay of correlations for predominantly expanding multimodal circle maps under moderate to large noise perturbations.  This class of systems consists of additive noise perturbations of deterministic maps whose derivative is very large in most of the phase space but in a small union of connected components where the derivative vanishes. Typically, in the absence of noise, the dynamics of these maps is trivial: almost surely all trajectories converge to a finite union of equilibria. This behaviour persists for small perturbations. However, as the noise strength exceeds a certain threshold, trajectories can escape the deterministic equilibria with high probability and spend more time in the expanding region, leading to positive Lyapunov exponent and chaotic behaviour \cite{Zeng}. This phenomenon is called {\em noise-induced chaos}.

The main novelty of our work is to develop a technique to build a random Young tower in a setting where the classical scheme, mentioned above, {\em does not work} for predominantly expanding circle endomorphism. Indeed, the systems we consider are noise-induced chaotic perturbations of deterministic maps that, in the absence of noise, have trivial dynamics, with their deterministic attractor being a union of sinks. This introduces a significant challenge: points can still become trapped with positive probability in the deterministic sink for arbitrarily long periods as the noise strength increases. Consequently, the essential finite-time mixing property required by the deterministic scheme fails.

To address this challenge, we introduce a novel concept: a sequence of stopping times, named \emph{Young times}, which captures the interaction between the frequent expansion guaranteed by hyperbolic times and the averaging effect of noise on the dynamics. This concept is specifically designed for {\em inherently noisy} systems, leveraging the chaotic effect generated by the noise to  describe the behaviour of orbit intervals.

Once the existence and main properties of Young times are established, the construction of the tower is completed through an adaptation of the classical scheme presented in \cite[Chapters 5 and 6]{Alves}. This construction retains several fundamental features of non-uniformly expanding maps and offers a versatile framework that could serve as a valuable tool for future developments.

{This work advances the study of predominantly expanding random circle endomorphisms. In \cite{Zeng,Blumenthal}, the authors established the positivity of the Lyapunov exponent and provided explicit lower bounds. Subsequently, \cite{Giuseppe} proved the existence of random horseshoes for these maps, which, in particular, implies positive quenched topological entropy. In this paper, we take this line of research further by showing that, as the expansion increases, such maps exhibit strong statistical quenched properties, namely quenched decay of correlations and the quenched almost sure invariance principle.  }

This paper is divided into five sections and two appendices. Section 2 presents the definition of predominantly expanding multimodal circle maps and states the main results  (Theorems \ref{TheoremB} and Corollaries \ref{decay} and \ref{QSAIP}). Section 3 focuses on defining Young times and proving their positive density.  In Section $4$, we construct a random Young tower for predominantly expanding multimodal circle maps. We also prove that the return function $R$ is measurable and show Theorem \ref{TheoremB}. In Section $5$, we prove Corollaries \ref{decay} and \ref{QSAIP}. Finally, in the final two appendices, we show Propositions \ref{LHiperbolicTijmens} and \ref{Thetainvmes}.

\section{Main results \label{main}}

\subsection{The model} \label{Model}

To ensure consistency with the predominant literature on multimodal circle maps, we adopt the definitions and the model presented in \cite{Blumenthal} throughout this paper. Specifically, we denote the unitary circle as $\mathbb S^1 = \mathbb R/\mathbb Z$, parametrized by the interval $[0,1)$. Also, we always assume that $\xi: \mathbb S^1\to \mathbb \mathbb \mathbb S^1$ is a \emph{non-degenerated $\mathcal C^2$ function}, i.e.
\begin{enumerate}
    \item[$(a)$] the critical set $\mathscr{C}:=\{x\in \mathbb S^1; \xi'(x) = 0\}$  has finite cardinality, and
    \item[$(b)$] $\{\xi''=0\}\cap \mathscr{C} = \emptyset.$
\end{enumerate}

Given $\alpha >0$ and $a\in [0,1)$ we define the family of maps
\begin{align}
    f = f_{\alpha,a} &:= \alpha \xi + a \ (\mathrm{mod }\ 1), \label{functionf}
\end{align}
where $(\mathrm{mod}\ 1)$ denotes the natural projection from $\R$ onto $\mathbb S^1.$ 

In the following, we induce a random dynamical system using the map \eqref{functionf}.  
Fix $\varepsilon > 0$ and define the space $\Omega = \Omega^{\varepsilon} := [-\varepsilon, \varepsilon]^{\mathbb{Z}},$ whose elements are two-sided sequences $\omega = (\ldots, \omega_{-1}, \omega_0, \omega_1, \ldots)$.  
Let $\nu^{\varepsilon}$ denote the uniform probability measure on $[-\varepsilon, \varepsilon]$, and define the product measure $\mathbb{P} = \mathbb{P}^{\varepsilon} := (\nu^{\varepsilon})^{\otimes \mathbb{Z}}$ 
on $\Omega$.  
We consider the $\sigma$-algebra $\mathcal{F} = \mathcal{B}([-\varepsilon, \varepsilon])^{\otimes \mathbb{Z}},$
where $\mathcal{B}([-\varepsilon, \varepsilon])$ denotes the Borel $\sigma$-algebra on the interval $[-\varepsilon, \varepsilon]$.

Given $\omega_0 \in [-\varepsilon,\varepsilon],$ we define $f_{\omega_0}(x) := f(x+\omega_0),$ where $x+\omega_0$ should be understood as $x + \omega_0\ (\mathrm{mod}\ 1).$ Now, we define the random dynamical system we aim to study. Setting $\alpha>0,$ $a\in \mathbb [0,1)$ and $\varepsilon >0$, for every $n\in \mathbb N$ we define the (random) iterations $ f_{\omega}^n :\mathbb S^1 \to \mathbb S^1$ as
\begin{align*}
f_{\omega}^n := f_{\omega_{ n-1}} \circ \ldots \circ f_{\omega_0}.
\end{align*}
Considering the shift map $\theta: \Omega \to \Omega,$ $\theta  ((\omega_i)_{i\in\mathbb Z} ) = (\omega_{i+1})_{i\in\mathbb Z},$ it is readily verified that for every $n,m \in\mathbb N$ 
$$f^{m+n}_{\omega}(x) = f^{n}_{\theta^m \omega} \circ f^m_{\omega}.$$
For every $\omega \in \Omega,$ we define  $\mathscr{C}_{\omega} := \left\{x\in\mathbb S^1;\  f'_\omega(x) = 0\right\}$ as the critical points of ${f_\omega = f_{\omega_0}.}$

Finally, we define the skew product of $\theta$ and $f_{\omega}$ as the map
\begin{align}
    \Theta: \Omega\times \mathbb S^1 &\to \Omega \times \mathbb S^1\label{skewproduct}\\
    (\omega, x)&\mapsto (\theta \omega, f_{\omega} (x) ).\nonumber
\end{align}
Observe that for every $n\in\mathbb N,$ $\Theta^{n}(\omega, x)= (\theta^n \omega, f^n_{\omega}(x)).$

{
The map \( (\omega, x, n) \in \Omega \times M \times \mathbb{N} \mapsto f_\omega^n(x) \in M \) defines a \emph{random dynamical system (RDS)} over the base dynamics \( (\Omega, \mathbb{P}, \theta) \). By slightly abusing notation, we refer to this cocycle structure as the \emph{random dynamical system $f^n_{\omega}$ generated by \( f_{\omega_0} \)}.

 During the paper we make use of the following notations:
\begin{itemize}
    \item[(i)] We denote by \( m \) the Lebesgue measure on \( \mathbb S^1 \). When integrating a function \( \varphi : \mathbb S^1 \to \mathbb{R} \) with respect to \( m \), we write
\[
\int_{\mathbb S^1} \varphi(x) \, dx := \int_{\mathbb S^1} \varphi(x) \, m(dx).
\]
\item[(ii)] Given \( \varphi \in L^\infty(\mathbb S^1, m) \), we define
$$\|\varphi\|_\infty = \inf \left\{ \sup_{x \in M} |\varphi(x)| \; ; \; m(M) = 1 \right\}\ \text{and } \|\varphi\|_{\mathrm{Lip}} = \|\varphi\|_\infty + \sup_{x \neq y} \frac{|\varphi(x) - \varphi(y)|}{d(x, y)},$$
where \( d \) is the arclength distance on \( \mathbb S^1 \); for simplicity, we write \( |x - y| := d(x, y) \).
\item[(iii)] Given a topological space \( X \), we denote its Borel \( \sigma \)-algebra by \( \mathcal{B}(X) \). For a family of probability measures \( \{\mu_y\}_{y \in Y} \) on \( X \) and a probability measure \( \mathbb P \) on \( Y \), if the map \( y \mapsto \mu_y(A) \) is measurable for each \( A \in \mathcal{B}(X) \), then we define the product measure \( \mu(dx, dy) := \mu_y(dx)\mathbb P(dy) \) on \( Y \times X \) by
\[
\mu(A \times B) = \int_A \mu_y(B) \, \mathbb P(dy), \quad \text{for all } A \in \mathcal{B}(Y), B \in \mathcal{B}(X).
\]
\end{itemize}}

\subsection{Main results}

In what follows, we state the main theorem of this paper. The definition of a \emph{random Young tower} is deferred to Section~\ref{randomtowers}.

{
\begin{mtheorem}\label{TheoremB}
 Let $a \in [0,1)$ and $\alpha,\varepsilon>0$. Consider the random dynamical system $f_\omega^n$ generated by ${f_{\omega_0}(x) = \alpha \xi (x+ \omega_0) + a \ (\mathrm{mod}\ 1)}$ on $\mathbb S^1,$ where $\omega_0 \in [-\varepsilon,\varepsilon],$ as defined in Section \ref{Model}. Then, for every $\alpha$ sufficiently large and for $\varepsilon>\alpha^{c-1}$, there exists an interval $\Delta\subset \mathbb S^1$, such that $f_\omega$ admits a random Young tower on $\Delta.$ Moreover, there exists a function $C \in L^2(\Omega, \mathbb P)$ and $\gamma > 0,$ such that the return function $R$, satisfies $ m [R_\omega > n] \leq  C(\omega) e^{-\gamma n},$ where $m$ is the Lebesgue measure on $\mathbb S^1.$
\end{mtheorem}}
Theorem~\ref{TheoremB} is proved in Section~\ref{proofB}.  
As a consequence of this theorem, we derive two corollaries that describe the quenched statistical behaviour of the random dynamical system $f_\omega^n$.
Corollary~\ref{decay} establishes the quenched decay of correlations for $f_\omega^n$, while Corollary~\ref{QSAIP} proves that $f_\omega^n$ satisfies the quenched almost sure invariance principle, under suitable conditions on the parameters $\alpha$ and $\varepsilon$.

{
\begin{mcorollary}\label{decay} Assume the same setting and notation as in Theorem~\ref{TheoremB}, and let $c \in (0,1)$. Then for $\alpha$ sufficiently large and for $\varepsilon > \alpha^{c -1}$,  there exists an invariant measure $\mu (\d \omega, \d x) = \mu_{\omega}(\d x) \mathbb P(\d \omega)$ for the skew product $\Theta,$ where $\mu_{\omega}$ is absolutely continuous with respect to the Lebesgue measure on $\mathbb S^1$ for $\mathbb P$-almost every $\omega\in \Omega.$  Moreover, there exist $C(\omega)\in L^2(\Omega)$ and $\gamma>0$ such that for every bounded measurable function $\varphi:\mathbb S^1\to \mathbb R$ and Lipschitz function $\psi:\mathbb S^1\to \mathbb R,$ 
$$ \left|\int_{\mathbb S^1} \varphi \circ f_{\omega}^n(x) \,  \psi(x) \d x - \int_{\mathbb S^1} \varphi(x)  \mu_{\theta^n \omega}(\d x) \int_{\mathbb S^1} \psi(x) \d x\right| \leq C(\omega) e^{-\gamma n}\|\varphi\|_{\infty} \|\psi\|_{\mathrm{Lip}},$$
and
$$ \left|\int_{\mathbb S^1} \varphi \circ f_{\theta^{-n}\omega}^n(x) \, \psi(x) \d x - \int_{\mathbb S^1} \varphi(x)  \mu_{\omega}(\d x) \int_{\mathbb S^1} \psi (x)\d x\right| \leq C(\omega) e^{-\gamma n} \|\varphi\|_{\infty} \|\psi\|_{\mathrm{Lip}}.$$
\end{mcorollary}}

{\begin{mcorollary}\label{QSAIP} 
Assume the setting and notation of Corollary~\ref{decay}, and fix $c\in(0,1)$.
There exists $\alpha_0\ge 1$ such that, for all $\alpha\ge \alpha_0$ and all $\varepsilon>\alpha^{c-1}$, the following holds.
Let $\varphi:\Omega\times \mathbb S^1\to\mathbb{R}$ satisfying
\begin{itemize}
    \item[$(a)$] $\omega \in \Omega \mapsto \|\varphi_\omega\|_{\mathrm{Lip}}\in \mathbb R$ belongs to $L^p(\Omega,\mathbb P)$ for some $p\in(1,\infty]$, where $\varphi_\omega(\cdot) := \varphi(\omega,\cdot)$;
    \item[$(b)$] $\int_M \varphi_\omega \d \mu_\omega =0$ for $\mathbb P$-a.e. $\omega\in \mathbb P$.
\end{itemize}
Given $n\in\mathbb N$ and $\omega\in \Omega$ define $S_n(\omega) = \sum_{i=0}^{n-1} \varphi_{\theta^i\omega}\circ f_\omega^i$ and
$\sigma_n^2(\omega) := \int  (S_n(\omega))^2 \d \mu_\omega.$ Then:
\begin{itemize}
    \item[(i)] There exists a constant $\sigma^2 \geq 0$ such that $\sigma_n^2(\omega)/n \xrightarrow[]{n\to\infty} \sigma^2$ for $\mathbb P$-a.e. $\omega\in \Omega$;
    \item[(ii)] If $\sigma^2 > 0$, the \emph{quenched almost sure invariance principle} holds.  
That is, there exists $\delta_0 \in (0, 1/4)$ such that, for $\mathbb{P}$-a.e. $\omega \in \Omega$, there exists a probability space $(\widehat{\mathbb{S}^1}, \widehat{\mu}_\omega)$ extending $(\mathbb{S}^1, \mu_\omega)$ and a sequence of measurable functions $\widehat{S}_n(\omega) : (\widehat{\mathbb{S}^1}, \widehat{\mu}_\omega) \to \mathbb{R}$
satisfying the following properties:
\begin{itemize}
    \item The processes $\{S_n(\omega)\}_{n \ge 1}$ and $\{\widehat{S}_n(\omega)\}_{n \ge 1}$ have the same distribution;
    \item There exists a standard Brownian motion $B^\omega$ on $(\widehat{\mathbb{S}^1}, \widehat{\mu}_\omega)$ such that
    \[
    \widehat{S}_n(\omega) 
    - B^\omega_{\sigma_n^2(\omega)} 
    = \mathcal{O}\big(n^{\frac{1}{4} + \delta_0}\big)
    \quad \widehat{\mu}_\omega\text{-a.s.}.
    \]
\end{itemize}
   \item[(iii)] If $\sigma^2 =0$ then $\varphi_\omega$ is a coboundary,  i.e. there exists a function $g\in L^\infty (\Omega\times\mathbb S^1,\mu)$ such that $\varphi = g\circ \Theta - g\ \mu\text{-a.s.}.$
\end{itemize}
\end{mcorollary}}
{Corollaries \ref{decay} and \ref{QSAIP} are proved in Section \ref{proofA}.}
\begin{remark}
The quenched almost sure invariance principle is a strong statistical property that implies several other statistical results, including the quenched central limit theorem, the quenched law of the iterated logarithm, and the quenched functional central limit theorem (see \cite[Remark 2.7]{Su}). {We remark that the choice of a uniform distribution for $\nu^\varepsilon$ on $[-\varepsilon,\varepsilon]$ is made solely for simplicity.  
Our results remain valid under the more general assumption that there exists a constant $C > 0$ such that
$$\frac{1}{C\,\varepsilon}\, \mathbbm{1}_{[-\varepsilon,\varepsilon]}
\ \leq\ 
\frac{\mathrm{d}\nu^\varepsilon}{\mathrm{d} m} 
\ \leq\ 
\frac{C}{\varepsilon}\, \mathbbm{1}_{[-\varepsilon,\varepsilon]}.$$}
\end{remark}

We note that the above results are sharp in the following sense: if $0 < \varepsilon \leq 1/\alpha$, the conclusions may fail, as the system $f_\omega^n$ can exhibit sinks whose basins of attraction have positive Lebesgue measure.  Moreover, it is worth emphasizing that Theorem~\ref{TheoremB}, as well as Corollaries~\ref{decay} and~\ref{QSAIP}, do not rely on any specific features of the function $\xi$ other than its non-degeneracy.  
This highlights a fundamental distinction between the classical theory of deterministic dynamical systems and its stochastic counterpart.

\subsection{Random Young towers}\label{randomtowers}
The groundwork for the results in this paper is the theory of random Young towers. In this subsection, we define the concept of random Young towers following closely the definitions provided in \cite{H, Baladi, Du}.

For simplicity, consider $(\Omega,\mathcal F,\mathbb P)$ and $\theta: \Omega \to \Omega$ as in Section \ref{Model}. Let $M$ be a compact {smooth manifold}  endowed with a {Borel} probability measure $\rho$, and $\{g_\omega: M \to M\}_{\omega \in \Omega}$ a family of maps depending only on the $0^{\mathrm{th} }$ coordinate of $\omega.$ We say that $g_\omega$ \emph{admits a random  tower on $\Delta \subset M$} if there exists a measurable return function $R:\Omega\times \Delta\to \mathbb N$, and for $\mathbb P$-almost every $\omega \in \Omega$ there exists a $\rho$-mod $0$ countable partition $\mathcal P(\omega) =\{P_j(\omega)\}_{j\in \mathbb N}$ of $\Delta,$ such that $R_\omega= R(\omega,\cdot)$ is constant on each $P_j(\omega)$ and
$$g^{R}_\omega(x) := {g_{\omega}^{R(\omega,x)}(x)} =  g_{\theta^{R_\omega(x) -1}\omega} \circ \ldots \circ g_{\theta \omega} \circ  g_\omega (x) \in \Delta\ \text{for }\mathbb P\otimes \rho \text{-almost every }(\omega,x)\in \Omega\times \Delta. $$

We define the $i$-\emph{th return time} as the following recursive equation:
$$
    R^i_\omega (x) = R^i (\omega,x) := \begin{cases}
    R(\omega,x)&, \ \text{if}\ i=1,\\
    R^{i-1}(\omega,x) + R\left(\theta^{R^{i-1}(\omega,x)} (\omega), g^{R^{i-1} (\omega,x)}_{\omega}(x)\right)&,\ \text{if }i>1.
    \end{cases}
$$

Assuming that $g_\omega$ admits a random tower on $\Delta$, we can define the random tower of $\mathbb P$-almost every $\omega \in \Omega$ as
\begin{align}
    \Delta_\omega = \left\{(x,\ell) \in \Delta\times \mathbb Z_{\geq 0}; \ x\in \bigcup_{j\in \mathbb N} P_j(\theta^{-\ell} \omega), \ \ell \in \mathbb Z_{\geq 0}, 0 \leq \ell \leq R_{\theta^{-\ell} \omega}(x) -1\right\},\label{Dw}
\end{align}
and the random tower map
\begin{align*}
    F_\omega : \Delta_\omega &\to \Delta_{\theta \omega}\\
    (x,\ell) &\mapsto
    \begin{cases}
    (x,\ell +1), \ \text{if}\ \ell +1 < R_{\theta^{-\ell} \omega} (x)\\
    (g^R_{\theta^{-\ell} \omega}(x),0), \ \text{if}\ \ell +1 = R_{\theta^{-\ell} \omega}(x).
    \end{cases}
\end{align*}
{More explicitly, we highlight that $F_\omega(x,\ell) = \left(g_{\theta^{-\ell}\omega}^{R(\theta^{-\ell}\omega,x))}(x),0\right)$ in the case that $R(\theta^{-\ell}\omega,x) = \ell +1$.}

We refer to $\Delta_{\omega,\ell} := \{(x,\ell) \in \Delta_\omega\} := \{(x,\ell); R_{\theta^{-\ell} \omega}(x) > \ell\}$ as the \emph{$\ell$-th level of the tower.} Observe that $\Delta_{\omega,0} = \Delta \times\{0\},$ for every $\omega\in\Omega$. Considering $$\mathbbold{\Delta} = \{(\omega , (x,\ell)); \omega\in \Omega \ \text{and}\ (x,\ell)\in \Delta_\omega\},$$
we denote $F: \tri \to \tri$ by $F(\omega,(x,\ell)) = (\theta \omega, F_\omega (x,\ell)).$ Using $\mathcal P(\omega)$, we induce a partition $\widehat {\mathcal P}(\omega)$ on $\Delta_\omega$ setting
$$\widehat{\mathcal P}(\omega) :=\left\{   P\times \{\ell\};  P \in \mathcal P(\theta^{-\ell} \omega),\ \left. R\right|_{P} \geq \ell +1 , \ell \in \mathbb Z_{\geq 0}\right\}.$$

Moreover, for a given $(x,\ell) \in \Delta_\omega,$ we denote $\widehat{R}^i_\omega$ as the $i$-th return time to the base level of $\tri$ i.e.
\begin{align}
    \widehat{R}^i_\omega (x,\ell) = R^i_{\theta^{-\ell} \omega}(x) - \ell.\label{hatr}
\end{align}
{Also, we denote $\widehat{R}_\omega(x,\ell) = \widehat{R}^{1}_\omega(x,\ell).$}
Furthermore, the reference measure $\rho$ on $\Delta$ can naturally be lifted to the set $\Delta_\omega$. We call the lifted measure $\rho_{\omega}$.

Now, we define the \emph{separation time} $s: \tri \times \tri \to \mathbb Z_{\geq 0} \cup \{\infty\}$ as 
$$s( ({\omega^1},z_1) , ({\omega^2},z_2)) := \begin{cases}
    0, \ \text{ if }\omega_1\neq \omega_2,\\
    \min\left\{n \geq 0; F_\omega^{\widehat R^n} (z_1)\ \text{and } F_\omega^{\widehat R^n} (z_2)\ \text{lie in a distinct element of }\widehat{\mathcal P}(\theta^{ \widehat R^n_\omega(x)}\omega) \right\}.
\end{cases}$$
When ${\omega^1} = {\omega^2} = \omega,$ we simply denote $s_\omega (z_1,z_2) := s((\omega,z_1), (\omega,z_2)).$

We say that $g_\omega$ \emph{admits a random Young tower on} $\Delta$  if $g_\omega$ admits a {random tower} on $\Delta$ such that the following properties hold:
\begin{enumerate}
    \item[(P1)] {\bf{Markov}}: For each $P\in \widehat{\mathcal P}(\omega)$ the map $F_\omega^{\widehat{ R}} : P \to \Delta_{\theta^{\widehat{R}_\omega}\omega, 0}=\Delta\times\{0\}$ is $\rho_\omega$-mod $0$ bijection, {observe that the function $\widehat{R}_{\omega}(\cdot)$ is constant} on $P$;
     \item[(P2)]{\bf{Bounded distortion}:} There exist $C_0>0$ and $0<\beta<1,$ such that for $\mathbb P$-almost every $\omega \in\Omega$ and $x,y \in P \in \widehat{\mathcal P} (\omega)$
     $$\log \frac{|(F_\omega^{{\widehat R}})' (x) |}{|(F_\omega^{{\widehat R}})' (y) |} \leq C_0 \beta^{s_\omega (F^{\widehat R}(\omega,x), F^{\widehat R}(\omega,y)) }. $$
     {The derivative \( (F^{\widehat{R}}_\omega)' \) is taken with respect to the natural product differentiable structure inherited from \( \Delta \times \mathbb{N} \)}.
     \item[(P3)]{\bf{Weak expansion}:} $\widehat{\mathcal P}(\omega)$ is a generating partition for $F_\omega$, i.e., the diameter of $\bigvee_{j=0}^{n} F_\omega^{-j}\widehat{\mathcal P}(\theta^j\omega)$ goes to $0$ as $n\to\infty.$
     \item[(P4)]{\bf{Aperiodicity}:} 
{There exist $N \in \mathbb{N}$, integers $t_1, t_2, \ldots, t_N \in \mathbb{N}$ with $\gcd\{t_1, t_2, \ldots, t_N\} = 1$, and constants $\varepsilon_1, \varepsilon_2, \ldots, \varepsilon_N > 0$ such that, for $\mathbb{P}$-almost every $\omega \in \Omega$ and each $i = 1, 2, \ldots, N$, we have
$
\rho\big\{ x \in \Delta : R_\omega(x) = t_i \big\} > \varepsilon_i.
$}

     \item[(P5)] {\bf{Non-degeneracy of the return time function}:} $R$ is a stopping time, i.e. if $R_{\omega}(x) = n,$ and $\upsilon \in \Omega$ and such that $\omega _i = \upsilon_i,$ for every $0\leq i < n,$ then $R_{\upsilon}(x) =n$, and
     $$ \int_{\Omega \times \Delta} R(\omega,x) \mathbb P \otimes \rho (\d \omega, \d x)<\infty.$$
\end{enumerate}

\section{Young times}

In \cite[Chapter 6]{Alves} using hyperbolic times, the author introduces a technique for constructing a Young tower applicable to non-uniformly expanding maps. However, this method does not readily extend to random systems. To address this issue, we propose a solution by defining the concept of Young times. In broad terms, Young times can be described as hyperbolic times such that $f_\omega^n$ presents bounded distortion estimates when returning to the reference set. 

This section formally introduces the notion of Young times (see Definition \ref{YoungTijmen}). To this end, we first establish two essential lemmas which are crucial to the definition of Young times and the construction of the random Young tower. Recall that we denote by $m$  the Lebesgue measure on $\mathbb S^1$.

\subsection{Technical lemmas}

We proceed to prove two technical lemmas.

\begin{lema}
Let $ a  \in [0,1)$ and consider $f$ as in \eqref{functionf}. 
 Then, given $K_0>0$ and $\gamma \in [0,1)$, there exists $\alpha_1>0,$ such that for every $\alpha >\alpha_1,$
 $$m\left\{x \in \mathbb S^1; \ |f'(x)| < K_0\alpha^{\gamma}\right\} < \frac{2 K_0 \# \mathscr{C}}{\min_{c\in\mathscr{C}}|\xi''(c)|} \frac{1}{\alpha^{1-\gamma}}. $$\label{LGiu}
\end{lema}

\begin{proof}

  Recall that $f = \alpha \xi + a\ (\mathrm{mod}\ 1)$. Let $\tilde \xi: \mathbb R\to \mathbb R$ be a lift of $\xi:\mathbb S^1 \to \mathbb R.$ Define ${\widetilde{\mathscr C}'_\xi = \{x\in [0,1); x \ (\mathrm{mod}\ 1) \in \mathscr{C}\}.}$ Observe that for every $x\in \mathbb R$ and $y\in \widetilde{\mathscr C}'_\xi$, 
 \begin{align}
   \tilde \xi(x) = \tilde \xi(y) + \frac{\tilde \xi''(y)}{2}(x-y)^2 + \smallO(|x-y|^2). \label{LGiu2}   
 \end{align}
Therefore,  $|\tilde \xi'(x)| > |\tilde \xi''(y)| |x-y| -  |\smallO(|x-y|)|.$ Defining $c_0 := \min\left\{\left|\tilde \xi''(x)\right|; x\in \widetilde{\mathscr C}'_\xi\right \} >0,$ we obtain that there exists $\delta>0$ such that, if  $x\in \mathbb S^1$ satisfies $\dist(x, \mathscr{C}) < \delta$, then   ${ |\xi'(x)| \geq c_0 \mathrm{dist}(x,\mathscr{C})/2}.$

Let $c_1 := \min\{|\tilde \xi'(x)| ; \ \dist(x,\mathscr{C})\geq \delta\} >0.$ Therefore, for every $\alpha > \alpha_1:= (K_0/c_1)^{1/(1-\gamma)},$ we obtain that
 $$m\left\{x\in \mathbb S^1; \ | f'(x)| <K_0\alpha^{\gamma}\right\} \leq \frac{2 K_0  \# \mathscr{C}}{ \min_{c\in\mathscr{C}}{|\xi''(c)|} } \frac{1}{\alpha^{1-\gamma}}.$$
 \end{proof}

\begin{lema}\label{CentralLema}
Fix $c\in (0,1)$ and $x_0\in \mathbb S^1$ so that $f'(x_0)\neq 0$. Let $\alpha, \varepsilon >0$ and $a\in [0,1)$. Consider $f_\omega^n$ as the random dynamical system defined in Section \ref{Model}. Then, for $\alpha$ sufficiently large and for $\varepsilon>\alpha^{c-1},$ there exists  $\widetilde \varepsilon >0$ and $\delta_0 >0$, such that for every $\delta_1<\widetilde{\varepsilon}$, and for every interval $I\subset \mathbb S^1$, $ \delta_1 \leq |I|,$ there exists $L = L(\delta_1),$ 
$\varepsilon_0 = \varepsilon_0(\delta_1)$ and $\iota = \iota(\delta_1) >0$, such that
    $$\mathbb{P}\left\{ \begin{array}{l}
    \mbox{$\exists$ $\ell \in \mathbb N,$ $0<\ell\leq L$ and an interval $J \subset I$ s.t., $f^\ell _\omega(J) =B_{2 \delta_0}(x_0)$} \\
    \mbox{and $\mathrm{dist}(f^i_\omega (J), \mathscr{C}_{\theta^i\omega} )> \varepsilon_0$ for every {$0\leq i<\ell$}.}
  \end{array}\right\} > \iota. $$\label{expansion}
\end{lema}

\begin{proof}

Take $\alpha_2 > 0$ such that for every $\alpha > \alpha_2$
\begin{align}
  \delta_0 :=\left((\# \mathscr{C} +1) + \frac{\# \mathscr{C}}{ \min_{c\in\mathscr{C}}{|\xi''(c)|} } \right)\frac{1}{\alpha^{1/2}} <\frac{1}{4}.\label{defd0}  
\end{align}
We divide the remaining proof into three steps.

\begin{step}{1}
We show that there exists $\alpha_3 > \alpha_2$ such that for every $\alpha > \alpha_3$ and any interval $I_0\subset \mathbb S^1$, satisfying $$m(I_0) = \frac{ 2k}{\alpha^{1-c}}, \ \text{for some }k\in\mathbb N.$$
We have that
\begin{itemize}
    \item[$(a)$] if $m(I_0) < 1/2,$ then there exists $J\subset I_0 \text{ and } m(f(J)) = 2 m(I_0),$
and
$$ \left.f'\right|_{J} > \alpha^{c/2}.$$
    \item[$(b)$] if $1/2\leq m(I_0) \leq 1,$ then there exists $J\subset I_0,$ such that $f(J)$ covers $\mathbb S^1$ twice, i.e. for every $x\in\mathbb S^1,$ there exists $x_1,x_2 \in J$ with $x_1 \neq x_2$ such that $f(x_1) = f(x_2) = x.$ Moreover, we have that $\left.f'\right|_{J} > \alpha^{c/2}.$
\end{itemize}
\end{step}

Applying Lemma \ref{LGiu} with $\gamma = 1-c/2,$ we can find $\alpha_3 > \alpha_2$ such that for every $\alpha > \alpha_3,$
$$ m\left( |f'| < \alpha^{c/2}\right) \leq \frac{2 \# \mathscr{C}}{ \min_{c\in\mathscr{C}}{|\xi''(c)|} } 
  \frac{1}{\alpha^{1-c/2}},$$
  and
$$\frac{1}{\#\mathscr{C}+1} \left( \alpha^{c/2} - \frac{ \# \mathscr{C}}{ k \min_{c\in\mathscr{C}}{|\xi''(c)|} } \right)> 4.$$

Let us prove $(a)$. Observe that there exists a connected subset $J_0\subset I_0\cap \{|f'|> \alpha^{c/2}\}$ such that
$$m(J_0) = \frac{1}{\#\mathscr{C} +1} \left(\frac{2 k }{\alpha^{1-c}} - \frac{2 \# \mathscr{C}}{ \min_{c\in\mathscr{C}}{|\xi''(c)|}} \frac{1}{\alpha^{1-c/2}} \right). $$

From the mean value theorem, we obtain
\begin{align}
    m(f (J_0) ) &\geq \frac{\alpha^{c /2}}{\#\mathscr{C} +1}\left(\frac{2 k }{\alpha^{1-c}} - \frac{2 \# \mathscr{C}}{ \min_{c\in\mathscr{C}}{|\xi''(c)|} } 
  \frac{1}{\alpha^{1-c/2}} \right)\nonumber \\
  &=\frac{2k}{\alpha^{1-c}} \frac{1}{\#\mathscr{C}+1} \left( \alpha^{c/2} - \frac{ \# \mathscr{C}}{ k \min_{c\in\mathscr{C}}{|\xi''(c)|} } \right)\nonumber  \\
  &\geq 4 m(I).  \label{xt}
\end{align}
Hence, one can find an interval $J \subset J_0 \subset  I_0\cap \{|f'|>\alpha^{c/2}\}$ which satisfies the desired properties. This proves $(a)$.

Now, we prove $(b).$ Consider $\tilde f: \mathbb R \to \mathbb R,$ as a lift of the function $f :\mathbb S^1\to \mathbb S^1$. Moreover, let $\widetilde I_0$ be a lift of $I_0.$ From the same computations provided in \eqref{xt}, we obtain that there exists $\widetilde{J}_0\subset \widetilde I_0,$  $$\widetilde{J}_0\subset \widetilde I_0 \cap \{f'\geq \alpha^{c/2}\},$$
such that
$$m (\widetilde{J}_0) \geq 4 m(I_0) > 2. $$
We conclude the proof of the claim by defining $J$ as the projection of $\widetilde{J}_0$ on $\mathbb S^1$. This concludes Step 1.

\begin{step}{2} Let $\alpha > \alpha_3,$ $I_0 = [i_0-\alpha^{c-1}, i_0+ \alpha^{c-1}]$ for some $i_0\in \mathbb S^1,$ and $n_0 = \min\{n; \ 2^n m(I_0) \geq 1\}.$ Then for every $n > n_0,$ there exists $J'(\omega_1,\ldots, \omega_n)\subset I_0,$ such that
\begin{enumerate}
    \item  $\displaystyle f_{\omega_n} \circ \ldots f_{\omega_1} \circ f (J'(\omega_1,\ldots ,\omega_n)) $ covers $\mathbb S^1$ twice;
    \item $( f_{\omega_{n}} \circ \ldots \circ f_{\omega_1}\circ f)'(x) > \alpha^{cn/2},$ for every $x\in J'(\omega_1,\ldots, \omega_n);$ and 
    \item $\displaystyle\left( \frac{m(I_0) }{\alpha \|f'\|_{\infty}}\right)^{n+1} \leq m(J'(\omega_1,\ldots,\omega_n )) \leq \left(\frac{m(I_0) }{\alpha^{c/2}}\right)^{n+1}.$
\end{enumerate}

Moreover, the intervals $J'(\omega_1,\ldots,\omega_n)$ can be chosen in a way that 
$$J'(\omega_1,\dots,\omega_n) \subset J'(\omega_1,\dots,\omega_k)\ \text{for every }1\leq k \leq n.  $$

\end{step}

We prove the result by induction on $n$ starting at $n_0$. If  $n = n_0$, since $f_{\omega_i} = f(x-\omega_i),$ the result is easily verified using Step 1. 

Assume that the result is valid for $n \geq n_{0}$, we show that it is also valid for $n+1$. Fix a path $(\omega_1,\ldots, \omega_{n},\omega_{n+1})$. From the induction hypothesis, we obtain that for every interval $I,$ such that $m(I_0) = 2/\alpha^{1-c}$, given $\omega_1,\ldots, \omega_{n},$ we obtain that there exists $J'(\omega_1,\ldots,\omega_n) \subset I_0$ such that
$$f_{\omega_{n}} \circ \ldots \circ f_{\omega_1} \circ f (J'(\omega_1,\ldots,\omega_n))\ \text{covers }\mathbb{S}^1\ \text{twice} $$
and
$$\left(f_{\omega_{n}} \circ \ldots \circ f_{\omega_1} \circ f\right)'(x) > \alpha^{\frac{n c }{2}} \ \text{for every }x\in J'(\omega_1,\ldots,\omega_n).$$

Therefore,
$$ \left(\frac{m(I_0)}{\alpha\|f'\|_{\infty}}\right)^{n +1}\leq m \left( J'(\omega_1,\ldots,\omega_{n})\right) \leq \left(\frac{m(I_0)}{\alpha^{ c/2}  }\right)^{n + 1}. $$

From Step 1 $(b)$, we obtain that there exists $$J'(\omega_1,\ldots, \omega_{n+1}) \subset J'(\omega_1,\ldots, \omega_n),$$
which satisfies the desired properties, proving Step 2.

\begin{step}{3}
We conclude the proof of the lemma.
\end{step}

Define $D_0 := B_{2\delta_0}(x_0)$. Given $n\geq n_0,$ we obtain from the above step that  $$f_{\omega_n} \circ f_{\omega_{n-1}}\circ \ldots \circ f_{\omega_1} \circ f (J'(\omega_1,\dots,\omega_n)) $$
covers $\mathbb S^1$ twice. Thus, we can find an interval $J(\omega_1,\ldots, \omega_n) \subset J'(\omega_1,\ldots,\omega_n)$ such that
\begin{itemize}
    \item $\displaystyle f_{\omega_{n}} \circ \ldots \circ f_{\omega_1} \circ f(J(\omega_1,\ldots, \omega_n)) = D_0;$
    
    \item $ f_{\omega_\ell}\circ \ldots f_{\omega_1}\circ f(J(\omega_1,\ldots, \omega_n)) \subset \left\{x \in \mathbb{S}^1; \ \d \left(f^\ell_{\omega}\right)(x-\omega_0) > \alpha^{\ell c/2}\right\},$ for every $1\leq  \ell \leq n;$ and

    \item $ \displaystyle \left(\frac{4 \delta_0}{\alpha\|\xi'\|_{\infty}}\right)^{n +1}\leq m \left( J(\omega_1,\ldots,\omega_{n})\right) \leq  \left(\frac{4 \delta_0}{\alpha^{ c/2}  }\right)^{n + 1}. $
\end{itemize}

Since $\left(f_{\omega_{\ell}} \circ \ldots \circ f_{\omega_1} \circ f\right)'(x) > \alpha^{\frac{\ell c }{2}}$ for every $x\in J(\omega_1,\ldots\omega_n)$ we obtain that there exists $\varepsilon_0 = \varepsilon_0(\delta_0),$ such that
\begin{align*}
   \dist\left(f_{\omega_{\ell}}\circ \ldots \circ f_{\omega_1} \circ f(J'(\omega_1,\ldots,\omega_n),\mathscr{C}_{\theta^\ell \omega}\right)> \varepsilon_0. 
\end{align*}

Let us fix  $L=L(\delta_1)$ big enough such that
$$\left(\frac{4 \delta_0}{\alpha^{ c/2}  }\right)^{L+1} \leq \frac{\delta_1}{2}\ \text{and }L\geq n_0. $$

For simplicity, let us denote
$$P:=  \mathbb{P}\left\{ \begin{array}{l}
    \mbox{$\exists$ $\ell \in \mathbb N,$ $0<\ell\leq L$ and an interval $J \subset I$ such that}, f^\ell _\omega(J) =B_{2 \delta_0}(x_0)\\
    \mbox{and $\mathrm{dist}(f^i_\omega (J), \mathscr{C}_{\theta^i\omega} )> \varepsilon_0$ for every $0<i<\ell$.}
  \end{array}\right\}.$$
From the construction of $P$, we obtain that
\begin{align*}
     P&\geq \mathbb P[\omega \in \Omega;\ I + \omega_0 \supset J(\omega_1,\ldots, \omega_L) ].
\end{align*}
Therefore, for any interval $I\subset \mathbb S^1,$ such that $m(I) \ge  \delta_1$ with middle point $i_0$  we obtain that if $i_0 +\omega_0\in J(\omega_1,\ldots, \omega_L),$ then ${I+\omega_0 \supset J(\omega_1,\ldots, \omega_L)}.$  Therefore
$$ P \geq  \frac{4\delta_0}{ (\alpha\|\xi'\|_\infty)^{L+1}}.$$ 

This concludes the proof of the lemma.
\end{proof}

\subsection{Young times} \label{RT}

In this subsection, we finally define Young times. 
Before, we recall the definition of hyperbolic times. The following definition is adapted from  \cite[Section 6.1.1]{Alves}. We start by fixing a real number $0 <b < 1/2.$
\begin{definicao}\label{hyperbolic times}
 Given $\sigma \in (0,1)$ and $r>0$, we say that $n$ is a $(\sigma^{2},r)$-hyperbolic time for $(\omega,x) \in \Omega \times \mathbb{S}^1$ if for every $0 \le k < n$, 
\begin{align*}
\prod_{i=k}^{n-1}\left|f'_{\theta^i\omega} (f_{\omega}^i(x))\right|^{-1} \leq \sigma^{2(n-k)},\ \text{and } \mathrm{dist}_{r}(f^{n-k}_\omega(x), \mathscr{C}_{\theta^{n-k}\omega}) \ge \sigma^{b {k}}, 
\end{align*}
where, given a set $A$, \ $$\ \dist_r(x,A) := \begin{cases}
\mathrm{dist}(x,A)&,\ \text{if} \ \dist(x,A)\leq r;\\
1&,\ \text{otherwise}.
\end{cases}$$

\end{definicao}

The result below states a consequence of $n$ being a $(\sigma^2,r)$-hyperbolic time for $(\omega,x)$.  The following result {can be found in \cite[Proposition 4.9]{Vil} (see also \cite[Proposition 6.7]{Alves}).}

\begin{proposicao}[{\cite[Proposition 4.9]{Vil}}]\label{defVn}
Given $\sigma ,r>0$, there exist $\tilde{\delta}=\tilde{\delta}(\sigma,r)$ and  $C= C(\sigma,r)>0,$ such that for every $(\omega,x)\in \Omega\times \mathbb S^1$ for which $n\in\mathbb N$ is a $(\sigma^2,r)$-hyperbolic time and for every $0<\delta<\tilde\delta$ the following items hold:
\begin{enumerate}\label{HipBalls}
    \item[(i)] There exists an open neighbourhood $V_n(\omega,x)$ of $x$, such that $f^n_\omega$ maps $V_n(\omega,x)$ diffeomorphically to $B_{\delta} (f_\omega^n(x)).$
    \item[(ii)] for every $z,y\in V_n(\omega, x)$
    \begin{enumerate}
        \item[$a.$] $|f_{\omega}^{n-k}(z) -f_{\omega}^{n-k}(y)| \leq \sigma^{k} |f^n_\omega(z) - f^n_\omega(y)|,$ for every $1\leq k \leq n$;
        \item[$b.$] $\displaystyle \log \frac{|(f^n_\omega)' (z)|}{|( f^n_\omega)' (y)|} \leq C  |f^n_\omega(z) - f^n_\omega(y)|$.
    \end{enumerate}
\end{enumerate}
\end{proposicao}

The proposition below gives us the existence of $L$-\emph{sparse hyperbolic times}. To not break the text flow, this proposition is proved in Appendix \ref{AppendixA}. The proof of this proposition closely follows the technique outlined in \cite{Giuseppe}.

\begin{proposicao}\label{LHiperbolicTijmens}
Let $L$ be a natural number and consider the family of random variables $\{\tau_i: \Omega\times \mathbb S^1 \to\mathbb N\cup \{\infty\}\}_{i\in\mathbb N}$ defined as $$ \tau_i(\omega,x) :=\begin{cases}
\min\{n; \ n\ \text{is a }(\sigma^2,r)\text{-hyperbolic time for }(\omega,x)\}&, \ if\ i=1, \\
\min\{n; \ n > L+\tau_{i-1}(\omega,x)\ \text{is a }(\sigma^2,r)\text{{-}hyperbolic time for }(\omega,x)\}&,\ \text{if }i>1.
\end{cases} $$

Given $c\in (0,1)$, for every $\alpha>0$  sufficient large, and for $\varepsilon>\alpha^{c-1},$ there exists $0<\sigma < 1$ and $r>0$, independent of $L$,  such that 
$$\mathbb P \otimes m  \left[\bigcap_{i\in\mathbb N}\{\tau_i < \infty\} \right] = 1, $$
and there exists $\gamma >0,$ such that
$$\liminf_{n\to\infty} \frac{1}{n}\#\left(\{1,2,\ldots,n\} \cap \{\tau_i(\omega,x)\}_{i=1}^{\infty}\right) \geq \frac{\gamma}{L+1}\ {\text{for $\mathbb P\otimes m$-a.e. $(\omega,x)\in \Omega\times M$}}.$$
Moreover, the convergence of
$$\sup_{x\in\mathbb S^1}\mathbb P \left[\#\left(\{1,\ldots,n\} \cap \{\tau_i(\omega,x)\}_{i=1}^{\infty}\right) \leq   \frac{\gamma}{2 (L+1)}  n\right] \xrightarrow{n\to\infty} 0,$$
holds exponentially fast.

\end{proposicao}

In the following, we finally define Young times. For every $\alpha, \varepsilon>0$ satisfying the conditions of Lemma \ref{expansion} and Proposition \ref{LHiperbolicTijmens}, let us associate  \begin{align}
    \delta_0, \delta_1, \varepsilon_0, \iota, L >0 \ \text{and }x_0\in \mathbb S^1,\label{constants}
\end{align} such that:
\begin{itemize}
    \item $\delta_0=\left((\# \mathscr{C} +1) + \frac{\# \mathscr{C}}{ \min_{c\in\mathscr{C}}{|\xi''(c)|} } \right)/\alpha^{1/2}$ as defined in \eqref{defd0};  
    \item $\delta_1$ small enough such for {$\delta = 9\delta_1<1/2$}  the conclusions of Proposition \ref{HipBalls} hold;
    \item $\varepsilon_0= \varepsilon_0(\delta_1)$ as it follows from Lemma \ref{expansion};
    \item $\xi'(x_0) \neq 0;$ and
    \item $L,\iota$ the conclusion of Lemma \ref{CentralLema} holds.
\end{itemize}

For a given interval $I$, to reduce notation and improve readability, we denote
\begin{align}
    \mathcal E(I) := \left\{ \omega\in\Omega\ ; \begin{array}{l}
    \mbox{$\exists$ $\ell \in \mathbb N,$ $0<\ell\leq L$ and an interval $J \subset I$ s.t., $f^\ell _\omega(J) =B_{2 \delta_0}(x_0)$} \\
    \mbox{and $\mathrm{dist}(f^i_\omega (J), \mathscr{C}_{\theta^i\omega} )> \varepsilon_0$ for every $0<i<\ell$.}
  \end{array}\right\}.\label{Ei}
  \end{align}

\begin{definicao}
Let $\{\tau_i\}_{i=1}^{\infty}$ be the family of $L$-sparse $(\sigma^2,r)$-hyperbolic times defined on Proposition \ref{LHiperbolicTijmens}. For every $n\in\mathbb N$, consider the set
\begin{align}
    Y_n(\omega, x) := \left\{i\in \{1,\ldots,n\}\cap \{\tau_i(\omega,x)\}_{i=1}^{\infty}; \ \theta^{i}\omega \in \mathcal E\left(B_{\delta_1}\left(f^i_\omega(x)\right)\right)\right\}.
\end{align} 
We say that $i$ is a $(\sigma^2,r)$-Young time for $(\omega,x)$ if there exists $n\in\mathbb N,$ such that $i\in Y_n(\omega,x)$.
\label{YoungTijmen}
\end{definicao}

The following theorem establishes that the Young times have positive density. We remark that for each $n\in \mathbb N$, $\omega\in\Omega$ and $x\in \mathbb S^1$, it follows that $Y_n(\omega,x)\subset \mathbb N$ and $\#Y_n(\omega,x) \leq \#Y_{n+1}(\omega,x).$ 

\begin{teorema}\label{youngprop}
Let $Y_n(\omega)$ be as in \eqref{YoungTijmen} and $c \in (0,1)$. Then, for sufficiently large $\alpha > 0$ and $\varepsilon > \alpha^{c-1}$, there exists $\theta_1 > 0$ such that for every $x \in \mathbb{S}^1$,  
$$\mathbb{P}\left[ \liminf_{n \to \infty} \frac{\# Y_n(\omega, x)}{n} > \theta_1 \right] = 1. $$
{Moreover, there exists $K_1,\kappa>0$ such that $\mathbb P[\#Y_n(\omega,x) \leq \theta_1 n] \leq K_1 e^{-\kappa n}$} for each $n\in\mathbb N$ and every $x\in\mathbb S^1$.
\end{teorema}

\begin{proof}
In order to simplify the notation let us denote $S_n(\omega,x) := \{1,\ldots, n\} \cap \{\tau_i(\omega,x)\}_{i=1}^{\infty}$. {From the Borel-Cantelli lemma, it suffices to show that there exists $\lambda > 0$ such that
$$
a_n := \mathbb{P}[\#Y_n(\omega, x) \leq \lambda n]
$$
decays exponentially fast as $n \to \infty$, uniformly in $x \in \mathbb{S}^1$. This will immediately imply both statements of the theorem.}

From Proposition \ref{LHiperbolicTijmens} and setting $\upsilon := \gamma/(2L +2)$ we obtain that
$$ \mathbb P[\#Y_n(\omega,x) \leq   \lambda n, \#S_n(\omega,x) \leq \upsilon n] \leq  \mathbb P[ \# S_n(\omega,x) \leq \upsilon n]  \xrightarrow{n\to\infty} 0,$$
exponentially fast and uniformly on $x\in\mathbb S^1$. Since
\begin{align*}
a_n := \mathbb P[\# Y_n(\omega,x) \leq \lambda  n, \# S_n(\omega,x) \leq \upsilon n] +  \mathbb P[\# Y_n(\omega,x) \leq \lambda  n , \# S_n(\omega,x) \geq \upsilon n],
\end{align*}
in to prove the theorem, it is enough to show that, for $\lambda >0$ small enough, the quantity 
$$ \mathbb P[\# Y_n(\omega,x) \leq \lambda  n, \# S_n(\omega,x) \geq \upsilon n]  \xrightarrow{n\to\infty} 0,$$
exponentially fast. We divide the remaining proof into three steps.

\begin{step}{1}
We show that 
\begin{align}
  &\mathbb P \left\{ \begin{array}{l}
   \# Y_n(\omega,x) \leq \lambda  n,\\
    \# S_n(\omega,x) \geq \upsilon n
  \end{array}\right\} =& \sum_{ s = \lceil \upsilon n \rceil}^{n} \sum_{j = s - \lfloor\lambda  n\rfloor}^{s} \sum_{0\leq i_1 \leq \ldots \leq i_j \leq s }\mathbb P \left\{ \begin{array}{l}
   \tau_{i_1}(\cdot,x) ,\ldots, \tau_{i_j}(\cdot,x) \not \in  Y_n(\cdot,x),\\
    \tau_{s}(\cdot ,x) \leq n < \tau_{s+1}(\cdot,x)
  \end{array}\right\}\label{giuseppe2}
    \end{align} 

\end{step}

Observe that
\begin{align*}
   \mathbb P[\# Y_n(\omega,x) \leq \lambda  n, \# S_n(\omega,x) \geq \upsilon n] &= \sum_{ s = \lceil \upsilon n \rceil}^{n} \mathbb P[\# Y_n(\omega,x) \leq \lambda  n, \# S_n(\omega,x) = s].
\end{align*}

Note that $\omega \in \{ \# Y_n(\cdot ,x) \leq \lambda  n , \# S_n(\cdot ,x) = s\}, $ if and only if 
there exists $ j  > s -\lambda  n$, such that there exists $j$ distinct $L$-sparse hyperbolic times $i_1,\ldots,i_j \in S_n(\omega,x) \setminus Y_n(\omega,x)$. Therefore
\begin{align}
   \mathbb P[\# Y_n(\omega,x) \leq \lambda  n, \# S_n(\omega,x) \geq \upsilon n] &= \sum_{ s = \lceil \upsilon n \rceil}^{n} \mathbb P[\exists\ j > s-\lambda  n; \exists\ i_1,\ldots \leq i_j \in S_n \setminus Y_n].\label{giuseppe1}
\end{align}

Observe that a hyperbolic time $i \in S_n \setminus  Y_n,$ if and only if  $i\in S_n(\omega,x)$ and $$\theta^{i}(\omega) \not \in \mathcal E(B_{\delta_1} (f_\omega^{i} (\omega,x) ) \ \text{if and only if } \omega \not \in  \theta^{-i} \left(\mathcal E(B_{\delta_1} (f_\omega^{i}(x))\right).$$
Moreover, $\# S_n(\omega,x) = s,$ if and only if $\tau_s\leq n < \tau_{s+1}.$

Combining the above observation with \eqref{giuseppe1} we obtain \eqref{giuseppe2}. This completes the proof of Step 1.

\begin{step}{2} \it{ We show that for every $j > s-\lambda  n$ and $0\leq i_1\leq  \ldots \leq i_j \leq s$, we have that
$$\mathbb P \left\{ \begin{array}{l}
   \tau_{i_1}(\cdot,x) ,\ldots, \tau_{i_j}(\cdot,x) \not \in  Y_n(\cdot,x),\\
    \tau_{s}(\cdot ,x) \leq n < \tau_{s+1}(\cdot,x)
  \end{array}\right\} \leq (1-\iota)^{ j }. $$}\label{claimgiu1}
\end{step}

Observe that 
\begin{align}
   & \mathbb P \left\{\omega \in \Omega\, ;\,  \begin{array}{l}
   \tau_{i_1}(\omega,x) ,\ldots, \tau_{i_j}(\omega,x) \not \in  Y_n(\omega,x),\nonumber\\
    \tau_{s}(\cdot ,x) \leq n < \tau_{s+1}(\cdot,x)
  \end{array}\right\} \\=& \mathbb P \left\{\omega \in \Omega\,;\, \begin{array}{l}
   \omega \in \bigcap_{k=1}^j \theta^{-\tau_{i_k}(\omega,x) }  \left(\mathcal E(B_{\delta_1} (f^{\tau_{i_k}(\omega,x)} (\omega,x)))\right)^{^\mathsf{c}} ,\\
    \tau_{s}(\omega ,x) \leq n < \tau_{s+1}(\omega ,x)
  \end{array} \right\}\nonumber\\
  \leq& \mathbb P \left\{\omega \in \Omega\, ;
   \, \bigcap_{k=1}^j \theta^{-\tau_{i_k}(\omega ,x) }  \left(\mathcal E(B_{\delta_1} (f^{\tau_{i_k}(\omega ,x)} (\omega ,x)))\right)^{^\mathsf{c}}  \right\}\label{giuseppe3}.
\end{align}

Moreover
\begin{align}
P &:=\mathbb P \left[\omega\in \Omega\, ;\,
   \bigcap_{k=1}^j \theta^{-\tau_{i_k}(\omega,x) }  \left(\mathcal E(B_{\delta_1} (f^{\tau_{i_k}(\omega,x)} (\omega,x)))\right)^{^\mathsf{c}}  \right]\nonumber \\ &= \mathbb E \left[\prod_{k=0}^{j-1} \mathbbm 1_{\theta^{-\tau_{i_k}(\cdot,x) }  \left(\mathcal E(B_{\delta_1} (f^{\tau_{i_k}(\cdot,x)} (\cdot,x)))\right)^{^\mathsf{c}} } \cdot\mathbbm 1_{\theta^{-\tau_{i_j}(\cdot,x) }  \left(\mathcal E(B_{\delta_1} (f^{\tau_{i_j}(\cdot,x)} (\cdot,x)))\right)^{^\mathsf{c}} }\right]\nonumber\\
   &=\mathbb E\left[ \mathbb E \left[\left.\prod_{k=0}^{j-1} \mathbbm 1_{\theta^{-\tau_{i_k}(\cdot,x) }  \left(\mathcal E(B_{\delta_1} (f^{\tau_{i_k}(\cdot,x)} (\cdot,x)))\right)^{^\mathsf{c}} } \cdot\mathbbm 1_{\theta^{-\tau_{i_j}(\cdot,x) }  \left(\mathcal E(B_{\delta_1} (f^{\tau_{i_j}(\cdot,x)} (\cdot,x)))\right)^{^\mathsf{c}} }\right\vert \mathcal F_{\tau_{i_j}(\cdot,x) } \right]\right]\label{giuseppe4},
\end{align}
where   $\mathcal F_{\tau_j(\cdot,x)}$ is the $\sigma$-algebra of $\tau_j(\cdot,x)$-past.

Observe that since $\{\tau_i\}_{i=0}^{\infty}$ are $L$-sparse hyperbolic times, from the definition of $\mathcal E(I)$ we obtain that
$$\prod_{k=0}^{j-1} \mathbbm 1_{\theta^{-\tau_{i_k}(\cdot,x) }  \left(\mathcal E(B_{\delta_1} (f^{\tau_{i_k}(\cdot,x)} (\cdot,x)))\right)^{^\mathsf{c}} }  \text{ is }\mathcal F_{\tau_j(\cdot,x)}\text{-measurable}.$$
Moreover, from the strong Markov property, we obtain that for every $\tilde \omega \in \Omega$
\begin{align*}
    \mathbb E\left[\mathbbm 1_{\theta^{-\tau_{i_j}(\cdot,x) }  \left(\mathcal E(B_{\delta_1} (f^{\tau_{i_j}(\cdot,x)} (\cdot,x)))\right)^{^\mathsf{c}} }\mid \mathcal F_{\tau_{i_j}(\cdot,x)}\right](\tilde \omega) &= \mathbb E\left[ \mathbbm 1_{ \left(\mathcal E(B_{\delta_1} (f^{\tau_{i_j}(\tilde \omega,x)} (\tilde \omega,x)))\right)^{^\mathsf{c}} }\right]\leq 1-\iota.
\end{align*}
The last inequality follows applying Lemma \ref{CentralLema} with $I =  m(B_{\delta_1} (f^{\tau_{i_j}(\tilde \omega,x)} (\tilde \omega,x)))$.

Combining, the above inequality with $\eqref{giuseppe4}$ if follows that
$$ P \leq (1-\iota) \mathbb P \left[\omega\in \Omega\, ;\,
   \bigcap_{k=1}^{j-1} \theta^{-\tau_{i_k}(\omega,x) }  \left(\mathcal E(B_{\delta_1} (f^{\tau_{i_k}(\omega,x)} (\omega,x)))\right)^{^\mathsf{c}}  \right], $$
repeating the same procedure $j-1$ times, we have that
$$ \mathbb P \left\{ \begin{array}{l}
   \tau_{i_1}(\cdot,x) ,\ldots, \tau_{i_j}(\cdot,x) \not \in  Y_n(\cdot,x),\\
    \tau_{s}(\cdot ,x) \leq n < \tau_{s+1}(\cdot,x)
  \end{array}\right\}\leq P \leq (1-\iota)^{ j }.$$

This proves Step 2. 

\begin{step}{3}
We show that for $\lambda$ sufficiently small, $\mathbb P[\# Y_n(\omega,x) \leq \lambda  n, \# S_n(\omega,x) \geq \upsilon n]  \xrightarrow{n\to\infty} 0$ exponentially fast.
\end{step}
From Step \ref{claimgiu1}, equation \eqref{giuseppe2} and assuming that $\lambda < \upsilon/2$ we obtain that

\begin{align}
  \mathbb P \left\{ \begin{array}{l}
   \# Y_n(\omega,x) \leq \lambda  n,\\
    \# S_n(\omega,x) \geq \upsilon n
  \end{array}\right\} \leq& \sum_{ s = \lceil \upsilon n \rceil}^{n} \sum_{j = s - \lfloor\lambda  n\rfloor}^{s} \sum_{0\leq i_1 \leq \ldots \leq i_j \leq s }(1-\iota)^j\nonumber\leq  \sum_{ s = \lceil \upsilon n \rceil}^{n} \sum_{j = s - \lfloor\lambda  n\rfloor}^{s} \binom{s}{j}(1-\iota)^j\nonumber\\
  &\leq  \sum_{ s = \lceil \upsilon n \rceil}^{n} \lambda  n \binom{s}{ s - \lfloor\lambda  n\rfloor}(1-\iota)^{( \upsilon   - \lambda) n}\nonumber\leq  \sum_{ s = \lceil \upsilon n \rceil}^{n} \lambda n \binom{s}{  \lfloor\lambda n\rfloor}(1-\iota)^{( \upsilon   - \lambda ) n}\nonumber\\
    &\leq  \lambda n^2 \binom{n}{  \lfloor\lambda n\rfloor}(1-\iota)^{( \upsilon   - \lambda ) n}\nonumber\leq  \lambda n^2 \left(\frac{e n}{\lfloor\lambda n\rfloor}\right)^{\lfloor\lambda n\rfloor}(1-\iota)^{( \upsilon   - \lambda ) n}\nonumber\\
    &\leq  \lambda n^2 \left(\frac{2e }{\lambda}\right)^{\lambda n}(1-\iota)^{( \upsilon   - \lambda ) n} \nonumber \leq  \lambda n^2 \left(\left(\frac{2e }{\lambda}\right)^{\lambda}(1-\iota)^{( \upsilon   - \lambda ) } \right)^{n}.\nonumber 
    \end{align} 
Since
$$\lim_{\lambda \to 0}  \left(\frac{2e }{\lambda}\right)^{\lambda}(1-\iota)^{( \upsilon   - \lambda )} = (1-\iota)^\upsilon < 1,$$
we can choose $\lambda$ small enough such that $\lambda< 2\upsilon$ and 
$\frac{1 + (1-\iota)^\upsilon }{2}<1.$  Therefore,
$$\mathbb P \left\{ \begin{array}{l}
   \# Y_n(\omega,x) \leq \lambda n,\\
    \# S_n(\omega,x) \geq \upsilon n
  \end{array}\right\} \leq \lambda n^2\left(\frac{1 + (1-\iota)^\upsilon }{2}\right)^n  \xrightarrow{n\to\infty} 0,  $$
exponentially fast. This proves Step 3. Therefore,  setting $\theta_1 = \lambda$, we complete the proof of the theorem.
\end{proof}

\begin{corolario}\label{giucor}
Let $c\in (0,1)$. For $\alpha>0$ sufficiently large and for $\varepsilon>\alpha^{c-1},$ 
$$ m\left\{x\in \mathbb S^1; \liminf_{n\to\infty} \frac{\#Y_n(\cdot,x)}{n} > \theta_1\right\}  =1 \ \mathbb P\text{-almost surely}. $$
\label{Gcor}
\end{corolario}
\begin{proof}
From Theorem \ref{youngprop} we obtain that for every $x\in \mathbb S^1$
$$\mathbb P \left[\liminf_{n\to\infty} \frac{\#Y_n(\omega,x)}{n} > \theta_1\right] = 1,$$
integrating both sides with respect to $m(\d x)$ 
$$\int_{\mathbb S^1}\mathbb P \left[\liminf_{n\to\infty} \frac{\#Y_n(\omega,x)}{n} > \theta_1\right] \mathrm {d} x = 1.$$

From Fubini's theorem
\begin{align*}
    1&=\int_{\mathbb S^1}\mathbb P \left[\liminf_{n\to\infty} \frac{\#Y_n(\omega,x)}{n} > \theta_1\right] \mathrm {d} x =\int_{\Omega} m\left[ \liminf_{n\to\infty} \frac{\#Y_n(\omega,x)}{n} > \theta_1 \right] \mathbb P [\mathrm {d}\omega] .
\end{align*}
Since,
$$ m\left[ \liminf_{n\to\infty} \frac{\#Y_n(\omega,x)}{n} > \theta_1 \right] \leq 1\ \text{for every }\omega\in \Omega,$$
we obtain that 
$$m\left[ \liminf_{n\to\infty} \frac{\#Y_n(\omega,x)}{n} > \theta_1 \right]  =1\ \text{for }\mathbb P\text{-almost every }\omega\in\Omega.$$

\end{proof}

Motivated by the above theorem, we define
\begin{align}
    Y(\omega,x) := \bigcup_{n\in \mathbb N} Y_n(\omega,x), \label{defY}
\end{align}
and for every fixed $n\in\mathbb N$ we denote
\begin{align}
    \widetilde{H}_n(\omega) := \{x\in \mathbb S^1 ; n\in Y(\omega,x)\}. \label{defHn}
\end{align}
From Corollary \ref{Gcor} we have that $\mathbb P$-almost every $\omega\in\Omega,$ we have that $m$-almost every $x\in \mathbb S^1$ lies in infinitely many $\widetilde{H}_n(\omega).$

\section{Random Young towers}\label{RYT}

In this section, we construct a random Young tower for $f_\omega^n$. The construction is based on the definition of the positive density of the Young times and the techniques presented in \cite[Chapter 5]{Alves}. Recall the definitions of the positive numbers $\delta_0,\delta_1, L$ and $x_0\in\mathbb S^1$ in \eqref{constants}.

From Proposition \ref{defVn} and the definition of $\delta_1$, if $n$ is a $(\sigma^2,r)$-Young time for $(\omega,x)$ we define $V_n(\omega,x)$ as the connected component of $f^{-n}_\omega \left( B_{9\delta_1}( f^n_\omega(x)\right)$ which contains $x$. In a similar fashion, we define the sets $W_n(\omega,x)\subset \widetilde{W}_n(\omega,x)\subset V_n(\omega,x)$ as,
$$ W_n(\omega,x) := C_x\left[f^{-n}_\omega \left( B_{\delta_1}\left( f^n_\omega(x)\right)\right)\right]\ \text{and}\  \widetilde{W}_n(\omega,x) := C_x\left[f^{-n}_\omega \left( B_{3 \delta_1}( f^n_\omega(x)\right)\right],$$
where $C_x[A]$ denotes the connected component of $A$ which contains $x$.

Suppose that $n$ is a $(\sigma^2,r)$-Young time, then $\theta^{n} \omega \in \mathcal E(B_{\delta_1}(f_\omega^i(x)))$. Therefore there exists $0\leq \ell \leq L$, and a set $U \subset f^n_\omega(W_n(\omega,x))$ such that $f^\ell_{\theta^n \omega}$ maps $U$ diffeomorphically to $B_{2 \delta_0}(x_0)$. We denote by $U'\subset U$ the set such that $f^{\ell}_{\theta^n \omega}$ maps $U'$ diffeomorphically to $B_{\delta_0}(x_0).$ Finally, we denote 
$\widetilde{w}_{n,\ell}(\omega,x)$ the subset of  $\widetilde{W}_n(\omega,x)$ that $f^n_\omega$ maps diffeomorphically to $U$ and we define \begin{align}
   \label{fakeI3} w_{n,\ell}(\omega,x)\text{ as the subset of }\widetilde{w}_{n,\ell}(\omega)\text{ such that} f^n_\omega\text{ maps diffeomorphically to }U'.
\end{align}
When the context is clear, we denote $\widetilde{w}_{n,\ell} (\omega,x)$ and ${w}_{n,\ell}(\omega,x)$, simply by $\widetilde{w}_{n,\ell}$ and ${w}_{n,\ell}$.

Defining $\Delta := B_{\delta_0}(x_0)$ we construct a partition of $\Delta$ using the above notation. In the following, we prove a lemma which allows us to naturally 
 adapt the construction presented {\cite[Chapter 5]{Alves}} to the random  context.

\begin{lema}\label{i1i2i3}
Let $c\in (0,1)$. For sufficiently $\alpha >0$ large  and for $\varepsilon >\alpha^{c -1}$. {There exist constants $C_0,C_1,\sigma>0$ and $\delta_1>0$ satisfying $9\delta_1 < 1/2$, such that }
\begin{enumerate}
    \item[$(I_1)$] For $\mathbb P$-almost every $\omega\in \Omega$, $\{\widetilde{H}_n(\omega)\}_{n\in\mathbb N}$  are compact sets (see \eqref{defHn}). Moreover, $m$-almost every point in $\Delta$ belongs to infinitely many $\widetilde{H}_n(\omega)$'s.
    \item[$(I_2)$] For every $x\in \widetilde{H}_n(\omega)$, $V_n(\omega,x)$ is diffeomorphically mapped in a disk of radius $9\delta_1$ around $f^n_\omega(x)$, and for every $y,z \in V_n(\omega,x)$
    \begin{itemize}
        \item  $|f_{ \omega}^{n-k}(z) -f_{ \omega}^{n-k}(y)| \leq \sigma^{k} |f^n_\omega(z) - f^n_\omega(y)|,$ for every $1\leq k \leq n$;
        \item  $\displaystyle \log \frac{|(f^n_\omega)' (z)|}{|(f^n_\omega)' (y)|} \leq {C_0}  |f^n_\omega(z) - f^n_\omega(y)|$.
        \end{itemize}
    
    \item[$(I_3)$] For every $\widetilde{w}_{n,\ell}(\omega,x)$ and $y,z\in f^n_{ \omega}(\widetilde{w}_{n,\ell}(\omega,x) ),$
    \begin{itemize}
        \item $\displaystyle \frac{1}{C_1} |f_{\theta^n \omega}^j({y}) - f_{\theta^n \omega}^j(z)| \leq |f^\ell_{\theta^n \omega} (y) - f^\ell_{\theta^n \omega} (z)| \leq C_1 |y-z|,$ for all $0\leq j \leq \ell$;
        \item $\displaystyle\log \frac{| (f^\ell_{\theta^n \omega})' (z)|}{|(f^\ell_{\theta^n \omega})'(y)|} \leq C_1 |f^\ell_{\theta^n \omega} (y) - f^\ell_{\theta^n\omega } (z)|$.
    \end{itemize}
    
\end{enumerate}
\end{lema}
\begin{proof}
Observe that $(I_1)$ follows directly from Corollary \ref{giucor}, and $(I_2)$ from Proposition \ref{defVn}.

We proceed to prove $(I_3)$. Since,  $f_\omega(x) = \alpha \xi_\omega(x+ \omega_0)+ a$, given $\varepsilon_0>0$ as in \eqref{Ei}, we can find a constant $K_0>0$ satisfying
$$\inf_{x\in S^1 \setminus B_{\varepsilon_0}(\mathscr{C}_{\omega} )}|d f_\omega(x)| > K_0,$$
where $B_{\varepsilon_0}(\mathscr C_\omega) = \{x\in \mathbb S^1; \mathrm{dist}(x,\mathscr{C}_{\omega}) <\varepsilon_0\}.$

Since $n$ is an $(\sigma^2,r)$-Young time, then $\theta^n\omega\in \mathcal E(B_{\delta_1}(f_\omega^i(x))).$ The above observation and the chain rule yields that there exists $\widetilde{K}_0$ such that
$$\frac{1}{\widetilde{K}_0} |f_{\theta^n \omega}^j(x) - f_{\theta^n \omega}^j(z)| \leq |f^\ell_{\theta^n \omega} (y) - f^\ell_{\theta^n \omega} (z)| \leq \widetilde{K}_0 |y-z|,$$ for all $0\leq j \leq \ell.$

To prove the second item of $(I_3)$, observe that
$$\displaystyle\log \frac{| (f^\ell_{\theta^n \omega})' (z)|}{|( f^\ell_{\theta^n \omega}) '(y)|} = \sum_{i=0}^{\ell-1} |\log|( f_{\theta^{\ell+i} \omega})' (f^{i}_{\theta^n\omega}z)  - \log|(f_{\theta^{\ell+i} \omega})' (f^{i}_{\theta^n\omega}y)| $$
using mean value theorem, the fact that the function $f$ is $\mathcal C^2$, and the first item of $(I_3)$, we find a constant $\widetilde K_1$ such that
$$ \log \frac{| (f^\ell_{\theta^n \omega)})'(z)|}{|(f^\ell_{\theta^n \omega})'(y)|} \le \widetilde K_1  |f^\ell_{\theta^n \omega} (y) - f^\ell_{\theta^n\omega } (z)|,$$
the proof is completed by defining $C_1 = \max\{\widetilde K_0, \widetilde K_1\}.$
\end{proof}

\subsection{Pathwise classical construction}\label{pathwiseconstr}

This subsection presents a pathwise tower construction for the random dynamical system $f_\omega^n$ on $\Delta$. It should be noted that the construction presented in this subsection does not result in a measurable return function $R$. However, in subsection \ref{measurabletower}, we slightly modify the construction presented in this section in order to obtain a measurable return function $R$.

The construction below is based on Lemma \ref{i1i2i3} and the method provided in {\cite[Chapter 5]{Alves}}. We closely follow the construction and the notation presented in {\cite[Chapter 5.2]{Alves}} to clarify and emphasise the connection between our work and the classical Young towers theory. 

\label{partition}
Fix a constant $c\in (0,1)$. Consider ${\alpha}> \alpha_4 = \alpha_4(c) >0$ large enough,  and $\varepsilon>\alpha^{c-1},$ such that the conclusions of Lemma \ref{i1i2i3} {hold} and
\begin{align}
    \frac{1}{\# \mathscr{C} +1}\left(2\delta_0 - \frac{2\# \mathscr{C}}{\min_{c\in \mathscr{C} }|\xi''(c)|}\frac{1}{\alpha^{1/2}} \right) \geq \frac{1}{\alpha^{1/2}}.\label{defa4}
\end{align}

For $\mathbb{P}$-almost every $\omega = (\omega_i)_{i\in\mathbb Z}$, we aim to construct {an} $m$-mod $0$  partition $\mathcal P(\omega)$ of $\Delta$ that can be used to build a Young Tower. To this end, we inductively define the families of sets $\{\mathcal P_n(\omega)\}_{n\in\mathbb N},$ $\{\Delta_n(\omega)\}_{n\in \mathbb N}$ and $\{S_n(\omega)\}_{n\in\mathbb N}$ depending only on $(\omega_0,\ldots, \omega_{n-1}),$ in a way that
\begin{itemize}
    \item $\mathcal P_n(\omega)$ is the family of elements of the partition of $\Delta$ constructed at step $n$;
    \item $\Delta_n(\omega)$ is the set of points which do not belong to  any element of the partition defined up to time $n$;
    \item $S_n(\omega)$ contains points in $H_n(\omega):= \widetilde{H}_n(\omega) \cap  \Delta$ (see \eqref{defHn}) not taken by elements of $\mathcal P_k(\omega)$'s constructed until time $n$.
\end{itemize}

In the following, we describe the construction of the above families of sets. Recall the constants $C_0, C_1$, $\delta_1$ and $\sigma$ defined in Lemma \ref{i1i2i3} and the constant $\delta_0$ from \eqref{defd0}. We start by defining the following auxiliary constants:
\begin{itemize}
    \item[(i)] $\delta_2 = (\delta_0 + 9\delta_1 C_1)/2$; 
    \item[(ii)] $1 < N_0\in\mathbb N$ large enough so that $C_1\sigma^{N_0} <1;$ and
    \item[(iii)] $N_1\in \mathbb N$ large enough so that
$$\delta_2 \sigma^{N_1} \leq \delta_0 \ \text{and }2\delta_1 \sigma^{N_1} \leq \frac{\delta_0(1-\sigma^{N_1})}{C_1}. $$
\end{itemize}
Also, given $w_{k,\ell}(\omega,x)$ we define for every $n>k,$
$$A_n(\omega, w_{k,\ell}) = \left\{y\in \widetilde{w}_{k,\ell};\ \mathrm{dist}(f^{k+\ell}_\omega(y), f_\omega^{k+\ell}( w_{k,\ell})) \leq \delta_0 \sigma^{n-k}\right\}. $$

From the definition of $A_n(\omega, w_{k,\ell})$ we obtain that $f^{k+\ell}_\omega (\widetilde{w}_{n,\ell})$ contains a neighbourhood of the boundary of $f^{k+\ell}_\omega(A_n(\omega))$ {whose length is at least}  $$2\delta_0 - \delta_0(1+\sigma^{n-k}) = \delta_0 (1-\sigma^{n-k}).$$

{We divide the construction into two parts: the first concerns the initial $N_0$ steps, and the second concerns the steps after $N_0$. Our approach follows \cite[Section 5.2.1]{Alves}, with a slight modification. In the construction presented in \cite[Section 5.2.1]{Alves}, nothing is done during the first $N_0 - 1$ steps, and the main construction is carried out at step $N_0$. In our case, we first construct a partition at the initial step and then do nothing until step $N_0+1$. Since the recurrent step is the same in both constructions, all subsequent arguments remain valid.}

   \begin{step}{1} If $n\in\mathbb N$  and $1\leq n \leq N_0.$ \label{step1}
\end{step}
For $n=1$ we define
$$\mathcal{P}_1(\omega) := \{ w_{(1,0)}(\omega, x(\omega)) \},$$
where ${w_{(1,0)}(\omega, x(\omega))}$ is an open subinterval of $\Delta$ such that
$$\left. f'_\omega \right|_{w_{(1,0)}(\omega, x(\omega))} \geq \alpha^{1/2}
\quad\text{and}\quad
m\big( \mathbb{S}^1 \setminus f_\omega(w_{(1,0)}(\omega, x(\omega))) \big) = 0.$$
Combining Lemma~\ref{LGiu} with \eqref{defa4}, we conclude that $\mathcal{P}_1(\omega)$ is non-empty.

We set
$$S_1(\omega, \Delta) :=
\big\{
f_\omega^{-1}\big(w_{(1,0)}(\omega, x(\omega))\big) \cap
B_{3\delta_1}\big(f_\omega(x(\omega))\big)
\big\},$$
and
$$S_1(\omega, \Delta^\mathsf{c}) :=
\big\{
x \in \Delta \; ;\;
\mathrm{dist}(x, \partial \Delta) < 2\delta_1 \sigma^{N_0}
\big\}.$$

Moreover, let
$$\Delta_1(\omega) := \Delta \setminus w_{(1,0)}(\omega, x(\omega)),$$
and define
$$S_1(\omega) := S_1(\omega, \Delta) \cup S_1(\omega, \Delta^\mathsf{c}).
$$

For definiteness, for $1 < n \leq N_0$ we set
$$\mathcal{P}_n(\omega) = \emptyset, \quad
S_n(\omega) = S_1(\omega), \quad
S_n(\omega, \Delta) = S_1(\omega, \Delta), \quad
S_n(\omega, \Delta^\mathsf{c}) = S_1(\omega, \Delta^\mathsf{c}), \quad
\Delta_n(\omega) = \Delta_1(\omega).$$

\begin{remark}
    {We note that the set $w_{(1,0)}(\omega, x(\omega))$ does not exactly coincide with the sets $w_{k,\ell}(\omega, y)$ defined in \eqref{fakeI3}. However, this slight abuse of notation will be convenient to ensure consistency when defining the return function $R:\Omega \times \mathbb{S}^1 \to \mathbb{N} \cup \{\infty\}$. Moreover, the set $w_{(1,0)}(\omega, x(\omega))$ is chosen in this way to ensure that the aperiodicity condition (P4) of the Young tower is satisfied later (see Step 1 in the proof of Theorem~\ref{TheoremB}).}
\end{remark}

\begin{step}{2}
 If $n\in\mathbb N$  and $ n > N_0.$
\end{step}
We proceed with the construction by induction. Give $n> N_0$, we assume that the families $\mathcal P_k(\omega)$, $S_k(\omega)$ and $\Delta_k(\omega)$ were constructed for every $1\leq k \leq n-1$. Then, we construct $\mathcal P_n(\omega),$ $S_n(\omega)$ and $\Delta_n.$

Let $F_n(\omega)$ be a finite subset of the compact $H_n(\omega)$ such that
$$H_n(\omega)\subset \bigcup_{x\in F_n(\omega)} W_n(\omega,x).$$

Consider $x_1, \ldots, x_{j_n} \in F_n(\omega)$ and, for every $1\leq i \leq j_n,$ a domain $w_{n,{\ell_i}}(\omega,x_i)\subset W_n(\omega, x_i)$ as in {\eqref{fakeI3}} for which $\mathcal P_n(\omega) = \{w_{n,\ell_1}(\omega,x_1) , \ldots, w_{n, \ell_{j_{n}}} (\omega,x_{\ell_n})\}$ is a maximal family of pairwise disjoint sets contained in $\Delta_{n-1}(\omega),$ such that for every $1\leq i \leq j_n$
\begin{align}
\omega_{n,\ell_i} \cap \left(\bigcup_{k=N_0}^{n-1} \bigcup_{w \in\mathcal P_k(\omega)} A_n(\omega, w)\right) = \emptyset.\label{partA}    
\end{align}
{To improve readability and reduce notation when the context is clear, we will simply write $w_{n,\ell_i}(\omega, x_i)$ as $w_{n,\ell_i}$.}

The sets in $\mathcal P_n(\omega)$ are the elements of the partition $\mathcal P(\omega)$ obtained in the $n$-th step of the construction. Set
$$\Delta_n(\omega)= \Delta \setminus \left(\bigcup_{k=N_0}^{n-1} \bigcup_{w \in\mathcal P_k(\omega)} w\right). $$

{Given $w_{k,\ell}(\omega,x) \in \mathcal P_k(\omega),$ for some $N_0 \leq k\leq n$, we define
$$S_n(\omega, w_{k,\ell}(\omega,x) ) = \widetilde W_{k}(\omega,x)\ \text{ if }n-k<N_1,  $$
and
\begin{align}
  S_n(\omega,  w_{k,\ell}) =\{y\in \widetilde{w}_{k,\ell} ; 0< \mathrm{dist}(f_\omega^{k+\ell}(y), f_{\omega}^{k+\ell}( w_{k,\ell}) \leq \delta_2 \sigma^{n-k}\}\ \text{if }n-k\geq N_1. \label{Sn}  
\end{align} 
We also set 
$$ S_n(\omega,\Delta)= \left(\bigcup_{k=1}^{n} \bigcup_{w \in \mathcal P_k(\omega)} S_n(\omega,w)\right)$$
and
$$S_n(\omega,\Delta^\mathsf{c}) =\left\{x\in \Delta; \, \mathrm{dist}(x,\partial \Delta^\mathsf{c}) < 2\delta_1 \sigma^n\right\}.$$}

Define
$$S_n(\omega)= S_n(\omega,\Delta^\mathsf{c})  \cup S_n(\omega,\Delta). $$

Finally, we write
\begin{align}
  \mathcal P(\omega) = \bigcup_{n\geq 1} \mathcal P_n(\omega). \label{defP}  
\end{align}

\begin{teorema}
{
Given $c \in (0,1)$, for every $\alpha > \alpha_4$ (see \eqref{defa4}) and $\varepsilon > \alpha^{c-1}$, the partition $\mathcal{P}(\omega)$ constructed above defines an $m$-mod 0 partition of $\Delta$.}

\end{teorema}

\begin{proof}

The proof follows directly from  \cite[Corollary 5.9]{Alves} applied to $\Delta_1(\omega)$.
\end{proof}

\subsection{Measurability of the return time \texorpdfstring{$R$}{R}\label{measurabletower}}

At this moment for $\mathbb P$-almost every $\omega\in \Omega,$ there exists {an} $m$-mod $0$ partition $\mathcal P(\omega)$ (see {\eqref{defP}}) of $\Delta = B_{\delta_0} (x_0).$ 
Following the construction presented in \cite[Chapter 5]{Alves}, we aim to define
\begin{align*}
    R:\Omega_0 \times \Delta &\to \mathbb N\\
    (\omega,x)&\to n+\ell\ \text{if }x\in\omega_{n,\ell}\in\mathcal P(\omega).
\end{align*}
However, it is unclear if $R$ is a measurable function.  To guarantee the measurability of the return function $R$, we slightly change the construction of $\mathcal P(\omega)$. 

We start by recalling the definition of a measurable family of sets.
\begin{definicao}
Given measurable spaces $(\Omega, \mathcal F )$ and $(X,\mathcal B(X))$, a family of sets ${\{A(\omega)\}_{\omega\in\Omega}\subset \mathcal B(X)}$ is called a \emph{measurable family of sets} if for every open set $U$, the set
    $$\{\omega\in \Omega; A(\omega)\cap U \neq \emptyset\}\text{ is measurable}. $$

When the context is clear we denote $\{A(\omega)\}_{\omega\in\Omega}$ simply by $A(\omega).$    
\end{definicao}

During this section, we make use of the following theorem.
\begin{teorema}[{\cite[Proposition 4.6]{vianaLya}}] \label{meas}Let $(\Omega,\mathcal F)$ and $(X,\mathcal B(X))$ be measurable spaces. Given a family of compact sets $\{A(\omega)\}_{\omega\in \Omega}\subset \mathcal B(X),$ the following statements are equivalent:
\begin{enumerate}
    \item $A(\omega)$ is a measurable family of sets; and
    \item the set $\mathrm{Graph}(A) =\{(\omega,x)\in \Omega\times X;\ x\in A(\omega)\}$ is $\mathcal F\otimes \mathcal B(X)$-measurable .
\end{enumerate}\label{eqv}
\end{teorema}

The following two results show that the sets used to construct the partitions in Section \ref{pathwiseconstr} form a family of measurable sets.

\begin{proposicao}
For every $n\in\mathbb N$, the set $A_n := \{(\omega,x)\in \Omega\times \mathbb S^1; V_n(\omega,x)\ \text{exists}\}$ is measurable. Moreover the families of sets $$\left\{H_n(\omega)\right\}_{\omega \in \Omega},\ \left\{\overline{V_n(\omega,x)}\right\}_{(\omega,x)\in \Omega}\ \text{and}\ \left\{\mathbb S^1 \setminus V_n(\omega,x)\right\}_{(\omega,x)\in\Omega\times M}$$are measurable. 
\end{proposicao}
\begin{proof}
Denote $\upsilon_i (\omega,x)$ as the $i$-th $(\sigma^2,r)$-Young time. Therefore
$$A_n = \bigcup_{i=0}^{\infty} \{\upsilon_i = n\} = \bigcup_{i=0}^{n} \{\upsilon_i = n\}, $$
which is measurable. Moreover, since $$H_n(\omega) =\Delta \cap  \bigcup_{i=0}^{n} \{\upsilon_i(\omega,\cdot) = n\},$$
we obtain that $H_n(\omega)$ is a measurable family of sets.

{

In the following, we show that the family $ \left\{\overline{V_n(\omega,x)}\right\}_{(\omega,x) \in \Omega \times \mathbb{S}^1}
$ is a measurable family of compact sets. For each $(\omega,x) \in A_n$, recall that $
f_\omega^n: V_n(\omega,x) \to B_{\delta_1}(f_\omega^n(x))$ is a diffeomorphism. From Proposition \ref{defVn}, we have $ V_n(\omega,x) \subset B_{9\delta_1 \sigma^n}(x). $ 

Since the map $(\omega,x) \mapsto f_\omega^n(x)$ from $\Omega\times\mathbb{S}^1$ to $\mathbb{S}^1$ is continuous, it follows that for each $n\in\mathbb{N}$,
$$
Z_n(\omega,x) := \overline{B_{9\delta_1\sigma^n}(x)} \cap (f_\omega^n)^{-1}\left(\overline{B_{\delta_1}(f_\omega^n(x))}\right)
$$
is a closed family of random compact sets. We also define the following families of random compact sets:
$$
Z_n^+(\omega,x) := Z_n(\omega,x) \cap [x,x+\delta_1],
\qquad
Z_n^-(\omega,x) := Z_n(\omega,x) \cap [x-\delta_1,x],
$$
where for $z,y \in \mathbb{S}^1$ we denote by $[z,y]$ the smallest arc in $\mathbb{S}^1$ connecting $z$ to $y$.

It follows that for each $(\omega,x) \in A_n$ we have
\begin{align*}
\overline{V_n(\omega,x)}
=& \left\{ \argmin_{z \in Z_n^{+}(\omega,x) \cap (f_\omega^n)^{-1}(\{p\})} |z-x| \,;\, p \in \overline{B_{9\delta_1}(f_\omega^n(x))},\ Z_n^{+}(\omega,x) \cap (f_\omega^n)^{-1}(\{p\}) \neq \emptyset \right\} \\
&\cup \left\{ \argmin_{z \in Z_n^{-}(\omega,x) \cap (f_\omega^n)^{-1}(\{p\})} |z-x| \,;\, p \in \overline{B_{9\delta_1}(f_\omega^n(x))},\ Z_n^{-}(\omega,x) \cap (f_\omega^n)^{-1}(\{p\}) \neq \emptyset \right\}.
\end{align*}

Observe that for each $(\omega,x)$, the set
$$
\{p \in \overline{B_{9\delta_1}(f_\omega^n(x))} \ \text{and} \ Z_n^{\pm}(\omega,x) \cap (f_\omega^n)^{-1}(\{p\}) \neq \emptyset\}
$$
is equal either to $[f_\omega^n(x),\, f_\omega^n(x) + \delta_1]$ or to $[f_\omega^n(x) - \delta_1,\, f_\omega^n(x)]$, depending on the sign of $(f_\omega^n)'(x)$.
Therefore, $\overline{V_n(\omega,x)}$ can be written as a finite union of inverse images of measurable functions in measurable compact sets. Since each $\overline{V_n(\omega,x)}$ is compact, we obtain a measurable family of compact sets.}

\end{proof}

\begin{corolario}\label{measurablesets}
For every $n\in \mathbb N$, the families of sets:
\begin{itemize}
    \item $\left\{\overline{W_n(\omega,x)}\right\}_{(\omega,x)\in \Omega\times \mathbb S^1}$ and $\{\mathbb S^1 \setminus W_n(\omega,x)\}_{(\omega,x)\in \Omega\times \mathbb S^1}$;
    \item $\left\{\overline{\widetilde W_n(\omega,x)}\right\}_{(\omega,x)\in \Omega\times \mathbb S^1}$ and $\{\mathbb S^1\setminus \widetilde{W}_n(\omega,x)\}_{\omega,x) \in \Omega\times \mathbb S^1};$
    \item $\left\{\overline{ w_{n,\ell} (\omega,x)}\right\}_{(\omega,x)\in \Omega\times \mathbb S^1}\ \text{and }\left\{\mathbb S^1 \setminus w_{n,\ell} (\omega,x)\right\}_{(\omega,x)\in \Omega\times \mathbb S^1}\text{for every }0\leq {\ell}\leq L;$ and
    \item $\left\{\overline{\widetilde w_{n,\ell} (\omega,x)}\right\}_{(\omega,x)\in \Omega\times \mathbb S^1}\ \text{and }\left\{\mathbb S^1 \setminus \widetilde w_{n,\ell} (\omega,x)\right\}_{(\omega,x)\in \Omega\times \mathbb S^1}\text{for every }0\leq {\ell}\leq L.$
\end{itemize}
are measurable. 
\end{corolario}
\begin{proof}
This follows from the definition of the above sets and the fact that $\left\{ \overline{V_n(\omega,x)}\right\}_{(\omega,x)\in \Omega\times \mathbb S^1}$  and $\left\{ \mathbb S^1 \setminus V_n(\omega,x)\right\}_{(\omega,x)\in \Omega\times \mathbb S^1}$  are measurable.
\end{proof}

The following three propositions show that the elements of $\mathcal P_n(\omega)$ can be chosen in a measurable way. {Observe that for each $\omega \in \Omega$, we let $\mathcal P(\omega)$ denote a measurable partition of $\Delta$. We often refer to the family $\{\mathcal P(\omega)\}_{\omega \in \Omega}$ as a random partition of $\Delta$. For brevity, we use the notation $\mathcal P(\omega)$ both to denote a single partition (for fixed $\omega$) and the family indexed by $\omega$, whenever no confusion arises.}

\begin{proposicao}
The partition $\mathcal P_{1}(\omega)$ can be taken measurably, i.e. there exist random variables $\delta,x:\Omega \to \Delta,$ such that
$$\mathcal P_{1}(\omega) = \{(x(\omega) -\delta(\omega), x(\omega) + \delta(\omega))\}_{\omega \in \Omega}. $$
\label{Giupar}
\end{proposicao}

\begin{proof}
Let us identify $\mathbb S^1$ with $[0,1)$, so that $\Delta$ is a connected interval when projected into $[0,1).$ By abuse of notation we denote $\Delta = [d_1,d_2]. $

Consider the random variable,
$$y(\omega) = \inf\left\{ y\in \Delta\cap [0,1) ; \ \inf_{q\in \mathbb Q \cap [0,1/\alpha^{1/2}] } |f'_\omega(y+q)| \geq \alpha^{1/2}\right\}. $$
Combining Lemma \ref{LGiu} with  \eqref{defa4} we obtain $y(\omega) \leq d_2 - 1/\alpha^{1/2}.$

Let us define the function $x:\Omega\to [0,1)$,\  $x(\omega) = y(\omega) + 1/(2 \alpha^{1/2}) \in \Delta$. Moreover, consider the function $\delta:\Omega \to [0,1),$ such that
$$\delta(\omega) = \inf\left\{0<r\leq 1/(2\alpha^{1/2}); \ m\left(\mathbb S^1 \setminus \left(f_\omega (-r+ x(\omega), x(\omega) +r )\right)\right) = 0\right\}. $$
Since $x,\delta:\Omega\to [0,1)$ are both random variables, the theorem is proved.

\end{proof}

\begin{proposicao}
Let $n\geq N_0$. In the construction of the partition $\mathcal P_{n}(\omega),$ the compact sets $\{F_n(\omega)\}_{\omega \in \Omega}$ can be taken measurably, i.e. the following properties hold:
\begin{enumerate}
    \item For every $n\in\mathbb N$ and $r\in\mathbb N_0$ the set
    $$M_{n,r} = \{\omega \in \Omega; \# F_n(\omega) =r\} $$
    is measurable, moreover
    \item   $\{F_n(\omega)\}_{\omega \in M_{n,r}}$ is a measurable family of compact sets.
\end{enumerate}
 
\label{GiuFn}

\end{proposicao}

\begin{proof}

Recall that we can identify $\mathbb S^1$ with the interval $[0,1)$. Consider the dyadic sets 
$$D_{i,j} := \left[\frac{i}{2^j}, \frac{i+1}{2^j}\right] \subset \mathbb S^1, j\in \mathbb N,\ 0\leq i \leq 2^{j} -1.$$

Since the family $\mathcal D = \{D_{i,j}\}_{j=0,0\leq i\leq 2^j - 1}^{\infty}$ is countable, we can reorder it in a way that $\mathcal D = \{D_\ell\}_{\ell \in \mathbb N}.$ Define the stopping times $i_{1,n} (\omega) :=  \min\{s\in \mathbb N;  H_n(\omega) \cap D_s \neq \emptyset\}.$ Observe that $i_{1,n}(\omega)$ is measurable since $H_n(\omega)$ is a measurable family of compact sets. Consider the random variable $X_{1,n}(\omega) := \inf\{ [0,1)\cap H_n(\omega) \cap D_{i_{1,n}(\omega) }(\omega)\}.$

For every $k \geq 1$, we recursively define stopping times
$$i_{k+1,n}(\omega) := \min\left\{s \in \mathbb N\,\left|\, \begin{array}{l} s \geq i_{k,n}(\omega),\   H_n(\omega) \cap D_s \neq \emptyset,\ \text{and}\\
\inf\{ [0,1)\cap H_n(\omega) \cap D_{s}\}\neq X_{j,n}(\omega),\ \forall \ 1\leq j\leq s
\end{array}\right.\right\},  $$
and the measurable random variable
$$X_{k+1,n}(\omega)  = \inf\{ [0,1)\cap H_n(\omega) \cap D_{i_{k+1,n}(\omega)}\}.$$

We claim that the function
$$k_n(\omega) = \inf\left\{r\in \mathbb N, H_n(\omega) \subset \bigcup_{i=1}^{r}  W_{n}(\omega,X_{i,n}(\omega) )\right\},$$
is  measurable. In fact, recall that $H_n(\omega)$ is compact $\mathbb P$-almost surely. Observe that for every $j \in \mathbb N$
$$H_n(\omega)\setminus \bigcup_{i=0}^{j-1} W_n(\omega,X_i(\omega)) =H_n(\omega)\cap \left(\bigcap_{i=0}^{j-1} \mathbb S^1\setminus W_n(\omega,X_i(\omega)) \right),$$
is a measurable family of compact sets. Finally, we obtain that
\begin{align*}
  \{k_n >j\} &= \left\{\omega;  H_n(\omega)\setminus \bigcup_{i=0}^{j-1} W_n(\omega,X_{i,n}(\omega)) \neq \emptyset \right\}\\
  &= \left\{\omega;  \left(H_n(\omega)\setminus \bigcup_{i=0}^{j-1} W_n(\omega,X_{i,n}(\omega)) \right) \cap  \mathbb S^1 \neq \emptyset \right\}.
\end{align*}
 is measurable.

From the above observation it follows that
$$F_n(\omega) := \{X_{1,n}(\omega),\ldots, X_{k_n(\omega),n} (\omega) \},$$
is a measurable family of compact sets.

\end{proof}

\begin{proposicao} \label{measurablepart}
For every $n\in \mathbb N$. The partition $\mathcal P_{n}(\omega)$ can be taken measurably, i.e. the following properties are fulfilled
\begin{enumerate}
    \item The set $\widetilde{M}_{n,r}= \{\omega\in \Omega ; \#\mathcal P_{n}(\omega) = r \} $ is measurable;
    \item For every $n,r\in\mathbb N$. There exists measurable function $\omega \in \widetilde{M}_{n,r}\mapsto (x_i(\omega), \ell_i(\omega))$ such that  $$\mathcal P_{n}(\omega) = \{ w_{n,\ell_1(\omega)}(\omega,x_1(\omega)) ,\ldots, w_{n,\ell_{r}(\omega)}(\omega,x_r(\omega))\},$$ and
    $$\left\{\overline{w_{n,\ell_i(\omega)}(\omega,x_i(\omega))} \right\}_{\omega \in \widetilde{M}_{n+1,r}} \text{ is measurable } \text{for every }1\leq i \leq r.$$
\end{enumerate}
\end{proposicao}
\begin{proof}
 We prove the result by induction on $n$. 
 By Proposition \ref{Giupar}, items $(1)$ and $(2)$ are true for every $n =1,\ldots, N_0-1.$ We show the result by induction on $n$. Suppose by strong induction that $(1)$ and $(2)$ is valid up to some $n\geq N_0-1$. We show that the result also holds for $n+1.$

Observe that by the induction hypothesis
$$ \omega \mapsto \widehat{A}(\omega) = \overline{ \bigcup_{k=1}^n \bigcup_{w\in \mathcal P_k(\omega)} A_{n+1}(\omega, w)},$$
is a measurable family of compact sets. 

Let us restrict ourselves to $ \omega \in M_{n+1,r} = \{\omega\in \Omega, \# F_{n+1}(\omega) = r\}$. Consider the set  $$G_{n+1,r} = \{(\omega,y_1,\ldots,y_r); \omega\in M_{n+1,r},\ y_i\in F_{n+1}(\omega)\}. $$
Since $\{F_n(\omega)\}_{\omega\in \Omega}$ is a measurable family of compact sets then, from Theorem \ref{eqv}, $G_{n+1,r}$ is a measurable set.

Observe that
\begin{align*}
    \widetilde{G}: G_{n+1,r}\times \{1,\ldots,L\}^r &\to \{1,2,\ldots, r\}\\*
    (\omega,(y_1,\ldots,y_r),(\ell_1,\ldots,\ell_r))&\mapsto \sum_{j=1}^{r} \left\lfloor \frac{ m\left(\widehat{A}(\omega) \cup \bigcup_{i=1}^{j} \overline{w_{n+1,\ell_i}(\omega, y_i)}\right) }{  \displaystyle m(\widehat{A}(\omega))+ \sum_{i=1}^{j} m\left(\overline{w_{n+1, \ell_i }(\omega,y_r) }\right)  }\right\rfloor
\end{align*}
is a measurable function. It is readily verified that $\widetilde{G}(\omega,(y_1,\ldots,y_n), (\ell_1,\ldots,\ell_r)) = s$ if and only if:
\begin{itemize}
    \item[$(i)$] $w_{n+1,\ell_i}(\omega, y_i) \cap w_{n+1,\ell_j}(\omega, y_j) = \emptyset, \ \text{for every}\ 1\leq i<j\leq s;$
    \item[$(ii)$] $ A(\omega)\cap w_{n+1,\ell_j}(\omega, y_j) = \emptyset, \ \text{for every}\ 1\leq i<j\leq s;$ and
\item[$(iii)$] $s$ is the largest natural number between $1$ and $m,$ such that $(i)$ and $(ii)$ holds. 

\end{itemize}

Define $G_{n+1,r}^r := \widetilde{G}^{-1}(\{r\}) \subset G_{n+1,r} \times \{1,\ldots, L\}^{r} \subset \Omega \times (\mathbb S^1)^r\times \{1,\ldots, L\}^{r}.$ Observe that $G_{n+1,r}^r$ is measurable. This implies that
$$K_{n+1,r,r}({\omega}):= \{(y, \ell ) \in F_{n+1}(\omega)^r\times \{1,\ldots, L\}^r; (\omega,y, \ell) \in G_{n+1,r}^r\}\subset F(\omega)^r\times\{1,\ldots,L\}^{r}, $$
is a measurable family of compact sets. Therefore $ L_{n+1, r, r} = \{\omega \in M_{n+1,r};\ K_{n+1,r,r}(\omega) = \emptyset\}$ is a measurable set.

From the above construction, given $\omega \in M_{n+1,r} \setminus L_{n+1, r, r}$, the function
$$(x,\ell)\in F_{n+1}(\omega)^r\times \{1,\ldots, L\}^{r} \mapsto \widetilde G(\omega, x,\ell), $$
attains its maximum if and only if $(x,\ell) \in K_{n+1,r,r}(\omega).$

Now we repeat the above process considering the set
$$G_{n+1,r}^{r-1} := \{(\omega,y); \omega\in L_{n+1,r,r}\ \text{and }y\in F_{n+1}(\omega)^r\}. $$
Observe that $ 0\leq \left.\widetilde G\right|_{G_{n+1,r}^{r-1}\times \{1,\ldots, L\}^{r-1}} \leq r-1.$ Therefore we can define the measurable set
$$G_{n+1,r,r-1} := \left(\left.\widetilde G\right|_{G_{n+1,r}^{r-1}\times \{1,\ldots, L\}^{r}}\right)^{-1}(\{r-1\}). $$
This implies that 
$$K_{n+1,r,r-1}(\omega) := \{(y,\ell) \in F_{n+1}(\omega)^r\times \{1,\ldots, L\}^{r}; (\omega,y,\ell) \in G_{n+1,r}^{r-1}\},$$
is a measurable family of compact sets. As before, we obtain that $$L_{n+1, r, r-1} = \{\omega \in L_{n+1, r, r}; K_{n+1,r,r-1}(\omega) = \emptyset\}, $$
is a measurable set. From the above construction, given $\omega \in L_{n+1, r, r} \setminus L_{n+1, r, r-1}$, the function
$$(x,\ell)\in F(\omega)^r\times \{1,\ldots, L\}^{r} \mapsto \widetilde G(\omega, x,\ell), $$
attains its maximum if and only if $(x,\ell) \in K_{n+1,r,r-1}(\omega).$

Finally, we define $G_{n+1,r}^{r-2} = \{(\omega,y);\, \omega \in L_{n+1, r, r-1}\ \text{and } y\in F_{n+1}(\omega)^r\}$ and repeat the above procedure $r-3$ times. In this way, we obtain the following $\mathbb P$-$\mathrm{mod}\ 0$ measurable decomposition of $M_{n+1,r},$
$$M_{n+1,r} = \bigsqcup_{i=1}^{r} \left(L_{n+1,r,i+1}\setminus L_{n+1,r,i}\right),$$
where $L_{n+1,r,r+1} = M_{n+1,r}$.

Consider the measurable family of compact sets
$$\widehat{G}_{n+1,r} (\omega) := \pi_{(\mathbb S^1)^i} \times \pi_{ \{1,\ldots, L\}^i}\left( K_{n+1, r , i}(\omega)\right) \ \text{if }\omega\in L_{n+1, r , i+1} \setminus L_{n+1,r,i}.$$
where $\pi_{(\mathbb S^1)^i} \times \pi_{ \{1,\ldots, L\}^i} ((y_1,\ldots,y_m),(\ell_1,\ldots,\ell_m)) = ((y_1,\ldots,y_i), (\ell_1,\ldots,\ell_i)).$ From \cite[Proposition 4.6]{vianaLya} there exists a measurable function map $$\sigma_{n+1,r,i}: L_{n+1,r,i+1}\setminus L_{n+1,r,i} \to (\mathbb S^1)^i\times \{1,\ldots, L\}^{i}$$ such that $\sigma_{n+1,r,i}(\omega) \in  \widehat{G}_{n+1,r} (\omega).$

Now, we show that $\mathcal P_{n+1}(\omega)$ satisfies $(1)$ and $(2)$. In order to see $(1)$, observe that 
$$\{\omega\in \Omega;\#\mathcal P_{n+1}(\omega) =r \} = \bigcup_{j=0}^{\infty} \left(L_{n+1, j ,r+1}\setminus L_{n+1, j ,r} \right),$$
which is measurable. Furthermore, item $(2)$ follows directly from the construction of the function $\sigma_{n,r, i}$ and Corollary \ref{measurablesets}.

\end{proof}

As a corollary of the above three propositions, we obtain that the partition $\mathcal P(\omega)$ of $\Delta$ can be modified such that the function $R$ defined at the beginning of this subsection is measurable.

\begin{corolario}\label{defR}
Let $\mathcal P(\omega)$ be a partition of the set $\Delta$ as constructed in subsection \ref{pathwiseconstr}, such that the conclusions presented in Proposition \ref{measurablepart} hold. Then, the return time function  \begin{align*}
    R:\Omega \times \Delta &\to \mathbb N\\
    (\omega,x)&\to n+\ell\ \text{if }x\in\omega_{n,\ell}\in\mathcal P(\omega),
\end{align*} is measurable. 
\end{corolario}
\begin{proof}
Recall from Proposition \ref{measurablepart} that
\begin{enumerate}
    \item the set $\widetilde{M}_{n+1,r}(\omega) = \{\omega ; \#\mathcal P_{n}(\omega) = r \} $ is measurable; and
    \item Given $\omega \in \widetilde{M}_{n,r}$. There exists measurable function $\omega \mapsto (x_i(\omega), \ell_i(\omega))$ such that  $$\mathcal P_{n}(\omega) = \{ w_{n,\ell_1(\omega)}(\omega,x_1(\omega)) ,\ldots, w_{n,\ell_{r}(\omega)}(\omega,x_r(\omega))\},$$ and
    $$\left\{\overline{w_{n,\ell_i(\omega)}(\omega,x_i(\omega))} \right\}_{\omega \in \widetilde{M}_{n,r}} \text{ is measurable } \text{for every }1\leq i \leq m.$$
\end{enumerate}

Observe that $\#\{\overline{w_{n,\ell_i(\omega)}(\omega,x_i(\omega))} \setminus {w_{n,\ell_i(\omega)}(\omega,x_i(\omega))} \} \leq 2.$ Therefore for every $k \in \mathbb N$
$$\{R=k\}= \bigcup_{n=1}^{k}\bigcup_{r=1}^\infty \bigcup_{i=1}^{r}\left\{(\omega,x); \omega \in \widetilde{M}_{n,r}, \ell_i(\omega) = k - n, x\in \overline{w_{n,\ell_i(\omega)}(\omega, x_{i}(\omega) )}  \right\}\ \left( \mathrm{mod }\ \mathbb P\otimes m\right).$$

Since $(\Omega\times \mathbb S^1,\mathbb P \otimes m)$ 
 is a complete probability space, we obtain that $R$ is a measurable function.
\end{proof}

\subsection{Random Young tower}
\label{proofB}
In this subsection, we prove Theorem \ref{TheoremB}. Namely, we show that for given $c\in(0,1)$, $\alpha>\alpha_4$ (see \ref{defa4}) and $\varepsilon> \alpha^{c-1}$, the random dynamical system $f_\omega$ admits a random Young tower.

{Using Proposition~\ref{measurablepart}, we fix, for $\mathbb{P}$-a.e.\ $\omega\in\Omega$, a partition $\mathcal{P}(\omega)$ of $\Delta$, chosen $m$-mod $0$.
Consequently, the return-time function $R_\omega:\Delta\to\mathbb{N}$ given by Corollary~\ref{defR} is measurable.
From now on, we work with these fixed choices; whenever we fix $\omega\in\Omega$, it is understood that $\omega$ lies in this full-measure set.}

We start proving the following lemma.
\begin{lema} \label{fr}
There are  constants $C>0$ and $0<\beta <1$ such that {for $\mathbb P$-a.e.} $\omega\in \Omega$ we have that for every $n\geq 1$ and $x,y \in w_{n,\ell} \in \mathcal P_n(\omega)$ the following holds:
\begin{enumerate}
    \item $|x-y|\leq \beta |f^{R}_\omega( x) - f^{R}_\omega(y) |$;
    \item $|f^j_\omega(x) - f^j_\omega(y)| \leq C|f^R_\omega(x) - f^R_\omega(y)|,\ \forall 0\leq j \leq R(\omega,x) = R(\omega,y);$ and
    \item $\displaystyle \log\frac{|(f_\omega^R)' (x)|}{|(f_\omega^R)' (y)|} \leq C |f^R_\omega (x) - f_\omega^R(y)|.$
\end{enumerate}
{Moreover, by increasing $\alpha$ and modifying the construction ${\mathcal P}(\omega)$, if necessary, $\beta \in (0,1)$ can be taken arbitrarily small.}
\label{panpilona}
\end{lema}

\begin{proof}

 First, observe that if $n = 1$, then $R(\omega, x) = R(\omega, y) = 1$ for each $x,y$ {lying in $\hat{w}(\omega,x(\omega))$ the unique element partition of $\mathcal P_1(\omega)$. Since $| f'_\omega(z)| \geq \alpha^{1/2} >1$ for every $z\in \hat w(\omega,x(\omega))$ and $\mathbb P$-a.e. $\omega \in \Omega$ we obtain the desired conclusion. Recall that $\mathcal{P}_n(\omega) = \emptyset$ for every $1 < n < N_0$.}

Now, suppose that $n\geq N_0$ and $w_{n,\ell} \in \mathcal P_n(\omega),$ we have that $w_{n,\ell} \subset V_n(\omega,z)$ for some $z \in H_n(\omega)$. Moreover $f^{n+\ell}_\omega$ maps $w_{n,\ell}$ diffeomorphically to $(x_0-\delta_0,x_0-\delta_0).$ From Lemma \ref{i1i2i3} $(I_2)$-$(I_3)$ we obtain that {for $x,y\in w_{n,\ell}$} we have that
$$|x-y| \leq \sigma^n |f^n_\omega(x) - f^n_\omega(y)| \leq C_1 \sigma^n | f^{n+\ell}_\omega(x) - f^{n+\ell}_\omega(y)|.$$
Setting $\beta = \max\{C_1\sigma^{N_0},\alpha^{-1/2}\} <1$, item $(1)$ is verified. {Observe that by increasing $N_0$ in Step 1 of Section~\ref{pathwiseconstr}, the constant $\beta$ can be made arbitrarily small, since the constant $C_1$ is uniform in $N_0$ (see Lemma~\ref{i1i2i3}).}

To prove item $(2)$, we separate it into two cases:
\begin{enumerate}
    \item[$(a)$] $n\leq j \leq n+\ell$; and
    \item[$(b)$] $0\leq j < n$.
\end{enumerate}
In case $(a),$ we obtain from Lemma \ref{i1i2i3}  ($I_3$) that
$$|f^j_\omega(x) - f_\omega^{j}(y)| \leq C_1 |f^{n+\ell}_\omega (x) - f^{n+\ell}_\omega (y)|.$$

On the other hand, if $(b)$ holds, from Lemma \ref{i1i2i3},  $(I_2)$ and $(I_3)$ we obtain that for every $x,y\in w_{n,\ell},$
$$|f^j_\omega(x) - f^j_\omega(y)| \leq \sigma^{n-j} |f_\omega^n(x) - f^n_\omega(y)| \leq C_1 \sigma^{n-j} |f^{n+\ell}_\omega (x) - f^{n+\ell}_\omega(y)|,$$
in both cases, we obtain the conclusion setting $C= C_1.$

The last one follows from Lemmas \ref{CentralLema} and \ref{i1i2i3} $(I_2)$-$(I_3)$ combined with
$$\displaystyle \log\frac{|(f_\omega^R)' (x)|}{|(f_\omega^R)' (y)|}\leq  \log\frac{|(f_{\theta^n \omega}^\ell)' (f^n_\omega x )|}{|(f_{\theta^n \omega}^\ell)' (f^n_\omega y )|} +\log\frac{|( f^n_{\omega})' (x )|}{|(f^n_{\omega})'(y)|} .$$

\end{proof}
{\begin{remark}
We note that, in this paper, the fact that $\beta \in (0,1)$ can be taken arbitrarily small is not used. This result is included since  this expansion factor estimation may still be relevant for particular applications or constructions.
\end{remark}}

Now, we prove Theorem \ref{TheoremB}.\\

\begin{theoremproofA}
    
    Fix $c\in (0,1)$. From the definitions provided in Section \ref{randomtowers}, we obtain that for every $\alpha > 
\alpha_4$ and $\varepsilon>\alpha^{c-1},$ the random dynamical system $f_\omega$ admits a random tower with partition $\mathcal P(\omega)$ and return time $R$. We apply the construction presented in Section \ref{randomtowers} setting $g_\omega= f_{\omega}$ and $\rho=m$.

We divide the proof into three steps.

\begin{step}{1} We show that $\mathrm{(P1)}$, $\mathrm{(P2)}$, $\mathrm{(P3)}$ and $\mathrm{(P4)}$ in Section \ref{randomtowers} hold.
    
\end{step}

(P1): Observe that (P1) follows directly from the definition of the sets $w_{n,\ell}(\omega,x)$ provided at the beginning of section \ref{pathwiseconstr}. 

(P2):  Let ${(\omega,z_1),(\omega,z_2)} \in \tri$. Assume that ${s_\omega(z_1,z_2) >0},$ otherwise, there is nothing to prove. {Recall the definition of $\widehat{R}^\ell$ for $\ell \in \mathbb N$ in \eqref{hatr}}. From Lemma \ref{panpilona} ($1$). 
\begin{align*}
    \left|F_\omega^{\widehat{ R}^1} (z_1) - F_\omega^{\widehat{ R}^1} (z_2) \right| &\leq \beta \left|F^{\widehat{ R}^2}_\omega(z_1) - F_\omega^{\widehat{ R}^2}(z_2)\right|\leq \ldots \leq \mathrm{diam}(\Delta) \beta^{s_\omega(z_1,z_2) -1}\\ 
    &=\mathrm{diam}(\Delta) \beta^{s_{\theta^{\widehat{ R}^1}\omega} \left(  F_\omega^{\widehat{ R}^1} (z_1), F_\omega^{\widehat{ R}^1} (z_2) \right)} = \mathrm{diam}(\Delta) \beta^{s\left(F^{\widehat{R}}(\omega,z_1) , F^{\widehat{R}}(\omega,z_2)\right)}.
\end{align*}
Using, Lemma \ref{panpilona} ($3$)
$$\log \frac{|(F_\omega^{{\widehat R}})'(z_1)|}{|(F_\omega^{{\widehat R}})'(z_2)|}\leq C |F^{\widehat R}_\omega(z_1) - F^{\widehat R}_\omega(z_2)| \leq \mathrm{diam}(\Delta)C  \beta^{s\left(F^{\widehat{R}}(\omega,z_1) , F^{\widehat{R}}(\omega,z_2)\right)}.$$

(P3): It is a standard consequence of  Lemma \ref{panpilona} (1) (for more details see \cite[Lemma 3.2]{Alves}).

(P4):  The construction of  $\mathcal P_1(\omega) \subset \mathcal P(\omega)$ carried out in Step \ref{step1} of Section \ref{pathwiseconstr} and Proposition \ref{Giupar} together imply  that for $\mathbb P$-almost every $\omega\in\Omega,$
$$m(x \in \Delta; R(\omega,x) =1\} \geq \frac{\mathrm{diam}(\Delta)}{\alpha \|\xi'\|_\infty}>0,   $$
{which implies that the aperiodicity condition $\mathrm{(P4)}$ is fulfilled with $N=1$.}

\begin{step}{2} We show that there exist $C(\omega)\in L^2(\Omega,\mathbb P)$ and $\gamma >0$ such that for $\mathbb P$-almost every $\omega\in\Omega$
$$m( R_\omega > n) \leq C(\omega) e^{-\gamma n},$$
where {we recall that }$R_\omega(\cdot) := R(\omega,\cdot).$
\end{step}

Let $\theta_1$ be as in Lemma \ref{CentralLema}, and define 
$$h_{\theta_1}(\omega, x ) := \min\left\{n\in\mathbb N; \frac{\#Y_i(\omega,x)}{i} \geq \theta_1,\ \text{for every }i\geq n\right\},$$
{where $Y_i(\omega,x)$ is was defined in Definition \ref{YoungTijmen}.}

Applying the same argumentation provided in  \cite[{Section} 5.3.2]{Alves} to each of the connected components of $\Delta_1(\omega)$ we are able to obtain that for $\mathbb P$-almost every $\omega\in \Omega$ and $n\geq {N_0},$ there exists a set $E_n(\omega)\subset \Delta,$ such that
\begin{align}
    \{R_\omega > n+ L \}\subset \{h_{\theta_1} > n\} \cup  E_n(\omega). \label{Rtail}
\end{align} 
{Moreover (see \cite[page 177-179]{Alves}), we can explicitly write
$$E_n(\omega):= \left\{x\in \Delta;\mathrm{dist}(x,\partial \Delta)\leq \delta_1\sigma^{\theta n/2}\right\}\cup X_{n,\left\lfloor \frac{\theta n}{2} \right\rfloor}\left(\omega\right),$$
where $\lfloor\cdot\rfloor$ denotes the floor function and given $\omega\in \Omega$, and $n,k\in\mathbb N$
\begin{align}
    X_{n,k}(\omega) :=\left\{x\in \Delta_n(\omega);\, \text{there exists }t_1<\ldots<t_k\leq n\ \text{such that }x\in \bigcap_{i=1}^k S_{t_i}(\omega,\Delta)\right\}. \label{Xnk}
\end{align}
From the computation developed in \cite[{Section} 5.3.2]{Alves} we observe that there exist constants $K,\kappa >0$, uniform in $\omega$, such that  $m(E_n(\omega)) \leq K_0 e^{-\kappa n}$ for every $n\in\mathbb N$. This is a consequence of $E_n(\omega)$ only depending on the definition of the satellite sets $S_i(\omega, \Delta)$. The decay estimates, in turn, rely solely on the geometric properties of the image of the sets $S_i(\omega,\Delta)$, which are uniform in $\omega$ due to the estimation of Lemma \ref{i1i2i3}.  We refer the interested reader to Remark \ref{Villarinho}, which clarifies why the constants in the above inequality are independent of $\omega$. }

    From \eqref{Rtail}, it is enough to estimate $m(\{h_{\theta_1}(\omega,\cdot) > n\}).$ Observe that
\begin{align*}
    \left\{(\omega,x) \in \Omega\times \mathbb S^1; h_{\theta_1}(\omega,x)> n\right\}
    &=\left\{(\omega,x)\in \Omega \times \mathbb S^1; \frac{\#Y_i(\omega,x)}{i} <\theta_1 \ \text{for some }i >n\right\}\\
    &\subset \bigcup_{i=n+1}^{\infty} \left\{(\omega,x)\in \Omega\times \mathbb S^1; \frac{\#Y_i(\omega,x)}{i}<\theta_1\right\}.
\end{align*}
From $\theta_1>0$ such that
$$\sup_{x\in\mathbb S^1}\mathbb P[\#Y_i(\cdot,x) < \theta_1 i]  < K_1 e^{-\kappa_1 i}.$$
Therefore, for every $x\in\mathbb S^1,$
$$ \mathbb P\left[ h_{\theta_1}(\cdot ,x)> n\right] \leq \sum_{i=n+1}^{\infty }K_{1}e^{-\kappa_1 i} \leq  \widetilde{K}_1 e^{-\kappa_1 n}\ \text{for some }{\widetilde{K}_1\geq 1}.$$

Observe that for every $\gamma \in (0,\kappa_1)$,
\begin{align}
\mathbb P\left[ m(h_{\theta_1}(\omega,\cdot)>j)> \widetilde{K}_1 e^{-\gamma j} \right] &\leq  \frac{1}{\widetilde K_1} e^{\gamma j } \mathbb E \left[m(h_{\theta_1}(\omega,\cdot)>j) \right]\nonumber\\
&\leq e^{ -(\kappa_1 -\gamma) j}.\label{gl}
\end{align}

From the Borel-Cantelli lemma, we obtain that the random variable,
$$n(\omega):=\max\left\{ j\in\mathbb N;\  m(h_{\theta_1}(\omega,\cdot)>j)> \widetilde{K}_1 e^{-\gamma j}\right\},$$
is $\mathbb P$-a.s. finite. Moreover, \eqref{gl} implies that
$$ \mathbb P[n(\omega) > k] \leq \sum_{i=k+1}^{\infty} \mathbb P\left[ m(h_{\theta_1}(\omega,\cdot)>i)> \widetilde{K}_1 e^{\gamma i } \right] \leq \frac{1}{1-e^{-(\kappa_1 - \gamma)}}e^{- (\kappa_1 - \gamma)k}.$$

Observe that for $\mathbb P$-almost every $\omega\in\Omega,$
\begin{align*}
    m(h_{\theta_1}(\omega,\cdot) > n) &\leq \widetilde K_1 e^{\gamma n(\omega) }e^{-\gamma n }.
\end{align*}
Setting $0< \gamma < \kappa_1/3$ and defining $C(\omega) := e^{\gamma n(\omega) }$,  we have that $C(\omega)\in L^2(\Omega)$ since
$$\int_{\Omega} C(\omega)^2\mathbb P[\d \omega] \leq  \sum_{i=0}^\infty \frac{e^{-(\kappa_1 - \gamma)i}}{1-e^{-(\kappa_1 - \gamma)}}e^{- (\kappa_1 - 3\gamma)i} <  \infty.$$
Finally, from \eqref{Rtail}  we obtain that
$$m(R_\omega > n +L) \leq C(\omega)e^{-\gamma n} + K_0 e^{-\kappa_0 n}, \ $$
which completes the proof of Step 2.
\begin{step}{3} We show $\mathrm{(P5)}$ in Section \ref{randomtowers} holds, which completes the proof of the theorem.
\end{step}
{Observe that, from the construction of the partition $\mathcal{P}(\omega)$ in Sections \ref{pathwiseconstr} and \ref{measurabletower}, as well as the definition of the return function $R$, it follows that $R$ is a stopping time. Indeed, the families of sets $\mathcal{P}_n(\omega)$, $\Delta_n(\omega)$, and $S_n(\omega)$ are constructed recursively, each step depending only on the past sequence $(\omega_0, \ldots, \omega_{n-1})$.}  Finally, from Step 2, we obtain that
     $$ \int_{\Omega \times \Delta} R(\omega,x) \mathbb P \otimes m (\d \omega, \d x)\leq \mathbb E \left[\sum_{i=0}^{\infty} n\, m(R_\omega \geq i)\right] \leq \mathbb E[C] \frac{e^{-\gamma}}{(1- e^{-\gamma})^2}<\infty.$$

This completes the proof of Step 3 and consequently proves Theorem \ref{TheoremB}.
\end{theoremproofA}

{
\begin{remark}[{Dependence of constants in \cite[Section 5.3.2]{Alves}}] \label{Villarinho}

The argumentation presented in \cite[Section 5.3.2]{Alves} shows that for $\mathbb P$-a.e. $\omega\in \Omega$, $$m(E_n(\omega))\leq m\left(\left\{x\in \Delta;\mathrm{dist}(x,\partial \Delta)\right\}\right)+m(X_{n,k}(\omega))\to 0$$ exponentially fast as $n\to \infty$ where the set $X_{n,k}(\omega)$ was defined in \eqref{Xnk}. The strategy of the proof consists in showing that the set $X_{n,k}(\omega)$ is contained in the union of suitable sets $Y_k(\omega)\cup Z_{n,k}(\omega)$. In particular, the argument shows that if $k=\lfloor \theta_1 n/2\rfloor$ then
$$
m\big(Y_k(\omega)\big)\ \le\ D_1\,\lambda_1^{\,k}
\ \text{and }\big(Z_{n,k}(\omega)\big)\ \le\ D_2\,\lambda_2^{\,k},
$$
with $0<\lambda_1,\lambda_2<1$ and constants $D_1,D_2>0$, possibly depending on $\omega\in \Omega$. 

We claim that the above constants can be taken uniformly on $\omega\in \Omega$ with probability one. To do so, we follow the procedure that leads to the construction of these constants and show that they depend only on the estimates of $(I_1)$-$(I_3)$ and on the length of suitably chosen uniformly small intervals.

Following the procedure in \cite[Section 5.3.2]{Alves}, we obtain constants $D_1>0$ and $\lambda_1\in(0,1)$ such that
$$
\sum_{\substack{n_1,\ldots,n_p \ge N_3 \\ n_1+\cdots+n_p = n}}
D_0^{\,p}\,\sigma^{\,n_1+\cdots+n_p}
\le D_1\,\lambda_1^{\,n},
$$
and $D_2>0$, $\lambda_2\in (0,1)$ with
$$
\frac{C_0^{\,2} C D\, D_2^{-1}}{1- D_2^{-1}} \le 1
\ \text{and }
\lambda_2 = \frac{D}{D+1},
$$
for suitable constants $D_0$, $C_0$, $C$, $D$, and a natural number $N_3$ (which may depend on $\omega$). We claim that these constants do not depend on $\omega$. Indeed, one can verify that
\begin{itemize}
\item $C_0$ is the distortion estimate coming from \cite[Lemma 5.2.]{Alves}   which depends on ($I_2$).
 \item $D_0$ is constructed in \cite[Lemma 5.14]{Alves} and is chosen such that $D_0\ge \max\{\widetilde C,\,\widetilde D\,C_0\}$, where $\widetilde C$ is the constant in \cite[Lemma 5.7]{Alves} which depends only on ($I_2$)-($I_3$), and $\widetilde D$ is an upper bound on the length of the union of at most two uniformly small intervals (since we are in dimension one).

  \item $N_3$ is chosen so that $\sigma + D_0\,\sigma^{N_3}<1$. 

  \item $D$ is defined using \cite[Lemma 5.7]{Alves} and depends only on ($I_2$)–($I_3$).

  \item $C$ is chosen in the proof of  \cite[Lemma 5.17]{Alves} as an upper bound on the length of the union of at most two uniformly small intervals (since we are in dimension one).

  \item  The proof also makes use of a natural number $N_4$ which is chosen to be large enough so that $D_2\,\sigma^{N_4}<1$  (see \cite[(5.72)]{Alves}).
\end{itemize}

\end{remark}}

\section{Quenched decay of correlations}
\label{proofA}
This section aims to prove Corollaries \ref{decay} and \ref{QSAIP}. For a fixed constant $c\in (0,1)$, we have established that for $\alpha >0$ sufficiently large and for $\varepsilon>\alpha^{c-1}$, the random dynamical system $f_\omega^n$ defined in Section \ref{Model}  admits a random Young tower with partition $\mathcal P(\omega)$ and return function $R$. Furthermore, we have shown that the probability of $\mathbb P\otimes m (R> n)$  converges exponentially fast to $0$ as $n$ approaches infinity. The forthcoming task involves linking the random Young tower properties established in Theorem \ref{TheoremB} with the original system.

{We begin by recalling a result from \cite{Du,Su} which states that, under the tower map $F : \tri \to \tri$, there exists a unique invariant measure absolutely continuous with respect to $m_\omega(\mathrm{d}x)\mathbb{P}(\mathrm{d}\omega)$ and having marginal $\mathbb{P}$. In our setting, this result takes the following form.

\begin{proposicao}[{\cite[Theorem 2.2.1]{Du},\cite[Lemma 3.1]{Su}}]\label{invmes}
Let $\alpha,\varepsilon >0$ such that Theorem \ref{TheoremB} holds. Then, the map $F:\tri \to \tri$ defined in Section \ref{randomtowers}  admits a unique invariant measure $\widetilde{\mu} = \widetilde{\mu}_\omega(\d x) \mathbb P(\d \omega)$ such that $\widetilde{\mu}_\omega \ll m_\omega$ and $\widetilde\mu_\omega ( \d x) /m_\omega(\d x) \in \mathcal F_\beta$ for $\mathbb P$-almost every $\omega\in \Omega$,  where $\beta$ is defined as in Proposition \ref{fr} and \begin{align}
    \mathcal F_\beta := \{\varphi:\mathbbold{\Delta} \to \mathbb R;\ \exists\ C \in L^{\infty}(\Omega,\mathbb P) \text{ s.t. } |\varphi_\omega(x) - \varphi_\omega(y)|\leq C(\omega) \beta^{s_\omega(x,y)}\ \text{for all }x,y\in \Delta_\omega\},\label{fbeta}
\end{align}
where $\varphi_\omega(\cdot)=\varphi(\omega,\cdot).$
\end{proposicao}

\begin{proof}
Let $\alpha, \varepsilon > 0$ be such that Theorem~\ref{randomtowers} holds.  
From Lemma~\ref{fr} we obtain  
$$1 < \frac{1}{\beta} \leq \inf \bigl\{ |(f_\omega^R)'(z)| : z \in \mathrm{Int}(P),\ P \in \mathcal{P}(\omega) \bigr\}
\quad\text{for }\mathbb{P}\text{-a.e. }\omega \in \Omega.$$
Moreover, the map $F : \tri \to \tri$ satisfies hypotheses (P1)–(P4) from Section~\ref{randomtowers}.  
Hence, by \cite[Theorem 2.2.1]{Du} (see also \cite[Lemma 3.1]{Su}), there exists a unique invariant measure   $\widetilde{\mu} = \widetilde{\mu}_\omega(\mathrm{d}x)\,\mathbb{P}(\mathrm{d}\omega)$
for $F$ on $\tri$, with $\widetilde{\mu}_\omega \ll m_\omega$ and  
$$\frac{\widetilde{\mu}_\omega(\mathrm{d}x)}{m_\omega(\mathrm{d}x)} \in \mathcal{F}_\beta
\quad\text{for }\mathbb{P}\text{-a.e. }\omega \in \Omega.$$

\end{proof}

}

Recall the skew product map $\Theta$ in \eqref{skewproduct}. Observe that the map 
\begin{align*}
   \pi: \tri &\to \Omega\times\mathbb S^1\\
   (\omega,(x,\ell))&\mapsto (\omega,f_{\theta^{-\ell}\omega }^\ell (x)),
\end{align*}
satisfies $\pi \circ F = \Theta \circ \pi$. Since $\pi$ is a semi-conjugacy between the maps $F:\tri \to \tri$ and ${\Theta: \Omega \times \mathbb S^1\to \Omega \times \mathbb S^1}$, when we consider the $F$-invariant measure $\widetilde{\mu}$ described in Proposition \ref{invmes}, we obtain that there exists a family of measures $\{\mu_\omega(\d x) \}_{\omega \in\Omega}$ on $\mathbb S^1$ such that
\begin{align}
    \mu(\d \omega,\d x) := (\pi_* \widetilde{\mu})(\d \omega,\d x)  = \mu_\omega(\d x) \mathbb P(\d \omega) \label{Thetainvmes}
\end{align}
is a $\Theta$-invariant measure. Observe that from the construction of $\pi$ we obtain that
$$\mu_\omega(\d x ) = \sum_{\ell =0}^{\infty} \widetilde \mu_{\theta^{-\ell}\omega}\left( ( f^{-\ell}_{\theta^{-\ell} \omega} (\d x) \cap \{R_{\theta^{-\ell} \omega }> \ell\}) \times\{0\}\right)\ \text{for }\mathbb P\text{-almost every }\omega\in\Omega.  $$
From the partition $\mathcal P(\omega)$ defined in Section \ref{pathwiseconstr} and Proposition \ref{Giupar}, and the regularity of the distribution of $\widetilde{\mu}$ over $m_\omega(\d x) \mathbb P(\d \omega),$ given in Proposition \ref{invmes}, it is readily verified that $\mu_\omega(\d x) \ll m(\d x)$ for $\mathbb P$-almost every $\omega\in\Omega.$ 

In the interest of conciseness, we provide a condensed version of \cite[Theorems 1.2.5 and 1.2.6.]{Du}. Despite being less general than the original theorems presented in \ref{Du}. The version below is enough to show Corollary \ref{decay}.

\begin{teorema}[{\cite[Theorem 1.2.5 and Theorem 1.2.6.]{Du}}]
Let $m$ be the Lebesgue measure on $\mathbb S^1,$ $\Omega$ be defined as in Section \ref{Model} and  $\{g_\omega: \mathbb S^1\to \mathbb S^1\}_{\omega\in \Omega}$ be a family of maps depending only on the $0^{\mathrm{th}}$ coordinate of $\omega$. Suppose that $g_\omega$ admits a random Young tower on $\Delta \subset \mathbb S^1$. Then there exists a unique exact absolutely continuous invariant probability measure $\nu(\d x, \d \omega) =\nu_\omega(\d x)\mathbb P(\d \omega)$ on $(\Omega\times \mathbb{S}^1, \mathcal F\otimes \mathcal B(\mathbb{S}^1) ),$ such that $\nu_\omega \ll  m$ for $\mathbb P$-almost every $\omega\in\Omega$.

Moreover, let $R:\Omega \times \Delta \to \mathbb N\cup \{\infty\}$ be the return time coming from the random Young tower structure of $g_\omega$ on $\Delta$, if there exists $K_0,\gamma_0,\kappa>0$ such that $$\mathbb P\otimes m(R>n) \leq K_0 e^{-\gamma_0 n}$$
and 
\begin{align}
    \inf\left\{|(g_\omega^R)'(x)|;\, x\in \mathrm{Int}(P)\ \text{and }P\in\mathcal P(\omega) \right\}>\kappa\ \text{for }\mathbb P\text{-a.e. }\omega\in \Omega.
\end{align}
Then, there exist $C(\omega)\in L^2(\Omega)$ and $\gamma>0$ such that for every bounded measurable function $\phi:\mathbb S^1\to \mathbb R$ and Lipschitz function $\psi:\mathbb S^1\to \mathbb R$ 
$$ \left|\int_{\mathbb S^1}\varphi \circ g_\omega^n \cdot \psi\, \d m - \int_{\mathbb S^1}\varphi \, \d \nu_{\theta^n \omega} \int_{\mathbb S^1} \psi\, \d m\right| \leq C(\omega) e^{-\gamma n}\|\varphi\|_{\infty} \|\psi\|_{\mathrm{Lip}}$$
and
$$ \left|\int_{\mathbb S^1}\varphi \circ g_{\theta^{-n}\omega}^n \cdot \psi\, \d m - \int_{\mathbb S^1}\varphi\, \d \nu_{\omega} \int_{\mathbb S^1} \psi\, \d m\right| \leq C(\omega) e^{-\gamma n} \|\varphi\|_{\infty} \|\psi\|_{\mathrm{Lip}}.$$

\label{Du}
\end{teorema}

Now we prove Corollaries \ref{decay} and \ref{QSAIP}.\\
\begin{theoremproofB}
    
    {Recall that  Lemma 4.12. implies that
    $$1< \frac{1}{\beta} \leq \inf \{|(f_\omega^R)'(z)|; z\in \mathrm{Int}(P)\ \text{and }P\in\mathcal P(\omega)\}\ \text{for }\mathbb P\text{-a.e. }\omega\in  \Omega. $$}
    Corollary \ref{decay} follows directly when combining Theorem \ref{TheoremB} and \ref{Du} when setting $g_\omega = f_\omega$ and $\nu= \mu$ (see \eqref{Thetainvmes}). {Corollary \ref{QSAIP} follows from Theorem \ref{TheoremB}  and \cite[Theorem 6.2]{Su} (see also \cite[Corollary 5.6]{liu2025quenchedinvarianceprinciplerate})}.
\end{theoremproofB}

\subsection*{Data availability statement} Data sharing not applicable to this article as no datasets were generated or analysed during the current study.

\subsection*{Conflict of interest statement} The authors have no competing interests to declare that are relevant to the content of this article.

\section*{Acknowledgments}

The authors wish to thank Alex Blumenthal for proposing the problem, and Jeroen S. W. Lamb,  Martin Rasmussen and Dmitry Turaev for their useful suggestions. MC’s research has been supported by an Imperial College President’s PhD scholarship and by the  Australian Research Council Laureate Fellowship (FL230100088). MC and GT are also supported by the EPSRC Centre for Doctoral Training in Mathematics of Random Systems: Analysis, Modelling and Simulation (EP/S023925/1).
\appendix

\section{Proof of Proposition \ref{LHiperbolicTijmens}}
\label{AppendixA}

In this appendix, we prove Proposition \ref{LHiperbolicTijmens}. We start by proving two technical lemmas.

\begin{lema}\label{largedeviazioni}
Let $a \in [0,1),$ and $\alpha,\varepsilon>0$. Consider the random dynamical system generated by ${f_{\omega}(x) = \alpha \xi (x+ \omega_0) + a \ (\mathrm{mod}\ 1)}$ on $\mathbb S^1,$ where $\omega_0 \in [-\varepsilon,\varepsilon],$ as defined in Section \ref{Model}. 
Then, for every $r$ small enough, there exists $\beta_0, r_0>0$ such that for every $x\in \mathbb S^1,$
\begin{align*}
  \mathbb{P}\left[ \sum_{i=0}^{n-1} -\log(\dist_r(f_\omega^i (x), \mathscr{C}_{\theta^i\omega})> n r_0 \right] \leq e^{-\beta_0 n} \ \text{for every }n\in\mathbb N.
\end{align*}
\end{lema}
\begin{proof}
 For every $n\in \mathbb Z_{\ge 0}$, consider $\mathcal F_n = \sigma(\pi_i; i \in \{0,1,\ldots, n-1\} ),$ where for every $i\in\mathbb Z$, we define $\pi_i: \Omega \to [-\varepsilon,\varepsilon]$ as $\pi_i( (\omega_n)_{n\in\mathbb Z}) = \omega_i.$

By the Markov inequality, given $1<\zeta< e$
\begin{align*}
  \mathbb{P}\left[ \sum_{i=0}^{n-1} -\log(\dist_r(f_\omega^i (x), \mathscr{C}_{\theta^i\omega})> n r_0 \right] \le \frac{\mathbb{E}\left[\prod\limits_{i=0}^{n-1} \zeta^{-\log\left(\mathrm{dist}_{r}\left(f^i_{\omega}(x),\mathscr{C}_{\theta^i\omega}\right)\right)}\right]}{\zeta^{r_0 n}}
\end{align*}
To simplify the notation, we define  write $b_i(\omega,x,r) := -\log\left(\mathrm{dist}_{r}\left(f^i_{\omega}(x),\mathscr{C}_{\theta^i\omega}\right)\right).$ Observe that 
\begin{equation}\label{induction}
\mathbb{E}\left[ \prod_{i\le n} \zeta^{b_i(\omega,x,\delta)}  \right] = \mathbb{E}\left[\prod_{i\le n-1}\zeta^{b_i(\omega,x,\delta)}\mathbb{E}\left[\zeta^{b_n(\omega,x,r)}|\mathcal{F}_{n-1}\right](\omega)\right]
\end{equation}
Since $\mathscr{C}_{\omega}=\mathscr{C}-\omega_0$ and the fact that $\mathscr{C}_{\theta^n\omega}$ depends explicitly on $\omega_n$, whilst $f^n_{\omega}(x)$ depends on $\omega_0,\dots,\omega_{n-1}$ we obtain  
\begin{align*}
\mathbb{E}\left[\zeta^{b_n(\omega,x,r)}|\mathcal{F}_{n-1}\right](\omega) &= \int_{[-\varepsilon,\varepsilon]} \zeta^{-\log\left(\mathrm{dist}_{r}(f^n_{\omega}(x)+\omega_n,\mathscr{C})\right)} \mathbb{P}(\d \omega_n)\\
&= \frac{1}{2\varepsilon} \int_{\mathbb{S}^1} \zeta^{-\log\left(\mathrm{dist}_{r}(y,\mathscr{C})\right)} 1_{B_{{\varepsilon}} (f^{n}_{\omega}(x))}(y)\d y \\
&=     \frac{1}{2\varepsilon}  \int_{\mathbb{S}^1} \mathbbm 1_{B_\varepsilon(f^n_{\omega}(x))}(y) \left(\frac{1}{\mathrm{dist}_{r}(y,\mathscr{C})}\right)^{\log(\zeta)} \d y  \\
&= \frac{1}{2\varepsilon}\int_{\mathbb{S}^1 \setminus B_r(\mathscr{C})}\mathbbm 1_{B_\varepsilon (f^n_{\omega}(x))}(y)\d y  + \frac{1}{2\varepsilon}\int_{B_r(\mathscr{C})}\mathbbm 1_{B_\varepsilon(f^{n}_{\omega}(x))}(y)\left(\frac{1}{\mathrm{dist}_{r}(y,\mathscr{C})}\right)^{\log(\zeta)}\d y 
\\
&\le  1 + \frac{2\# \mathscr{C}}{\varepsilon } \int_{(0,r)}\left(\frac{1}{x}\right)^{\log(\zeta)} \d x \\
&=   1 + \frac{2\# \mathscr{C}}{\varepsilon }\frac{r^{1-\log(\zeta)}}{1-\log(\zeta)},
\end{align*}
Substituting this estimate into \eqref{induction}, we achieve by induction on $n$ that
\begin{equation}\label{bondo3}
\mathbb{E}\left[ \prod_{i\le n} \zeta^{b_i(\omega,x,\delta)}  \right] \le \left(\frac{1}{\zeta^{r_0}}\left(1 + \frac{2\#\mathscr{C}}{\varepsilon }\frac{r^{1-\log(\zeta)}}{1-\log(\zeta)} \right)\right)^n.
\end{equation}

Thus, for $r$ small enough such that we obtain
$$\beta_0 := -\log\left( \frac{1}{\zeta^{r_0}}\left(1 + \frac{2\#\mathscr{C}}{\varepsilon }\frac{r^{1-\log(\zeta)}}{1-\log(\zeta)} \right) \right) > 0,$$
 which, from  \eqref{bondo3},  concludes the proof.
\end{proof}

\begin{lema}\label{ledev}

Let $a \in [0,1)$ and $\alpha,\varepsilon > 0$. Consider the random dynamical system generated by $f_{\omega}(x) = \alpha \xi (x+ \omega_0) + a \pmod{1}$  
on $\mathbb{S}^1$, where $\omega_0 \in [-\varepsilon,\varepsilon]$, as defined in Section~\ref{Model}.  

Let 
$C = \{x \in \mathbb{S}^1 : |f'(x)| < 1\}.$ 
Given $R > 0$, define  $
G = G(R) = \{x \in \mathbb{S}^1 : |f'(x)| > R\}.
$
Assume that there exist constants $R > 1$, $h \in (0,1)$, and $\varepsilon > (1-m(G))/h$ such that  
$$
Z(h) := \log(R)(1-h) + \frac{1}{2\varepsilon}\int_{C} \log|f'(z)|\,dz > 0.
$$  
Then, there exist constants $D, \beta_1, \sigma^2 > 0$ such that  
\begin{equation}\label{crazy2}
\sup_{x \in \mathbb{S}^1} 
\mathbb{P}\!\left( \sum_{i=0}^{n-1}\log\left|f'_{\theta^i\omega}(x)\right| < \sigma^2 n \right) 
\le D e^{-\beta_1 n}.
\end{equation} \label{new}

\end{lema}

\begin{proof}
We divide the proof into four steps. The proof is long and technical. In Step 1, we provide helpful definitions for the proof. Steps 2 through 4 involve estimating the measure of certain events, which subsequently leads to the desired result.
\begin{step}{1}
Initial preparation for the proof.
\end{step}

In order to simplify the notation, we denote $V:=-\int_C\log(|f'(x)|)\d x>0$.
Let $\beta_2>0$ such that
\begin{equation}\label{zetauno}
Z_1(h) := \log(R)(1-h)-\frac{\beta_2+1}{2\varepsilon}V.
\end{equation}
Take sufficiently large $\ell \in \mathbb{N}$ and consider the sets
\begin{equation}\label{defPk}
P(k):= \left\{x \in \mathbb{S}^1 \mid R^{-\frac{k}{\ell}}>|f'(x)|>R^{-\frac{(k+1)}{\ell}} \right\}.
\end{equation}
Observe that $\{ P(k)\}_{k \ge -\ell}$ is a partition of $\mathbb{S}^1 \setminus G$ whose diameter is small if $\ell$ is large enough and the sets $\{P(k)\}_{k \ge 0}$ are a partition of $C$  with the same properties. Then, for any sequence $\{x_k\}_{k \ge 0}$ such that $x_k \in P(k)$ we have that
\begin{align*}
 \left|\sum_{k \ge 0} \log|f'(x_k)|m(P(k)) +V\right| \xrightarrow{\ell \to\infty} 0.
\end{align*}
Furthermore, as a consequence of the definition of $P(k)$ in \eqref{defPk}, we have that  
\begin{align*}
\left|\sum_{k \ge 0} \log|f'(x_k)|m(P(k))- \sum_{k \ge 0} \log\left(R^{-\frac{(k+1)}{\ell}}\right)m(P(k)) \right| \le \log\left(R^{-\frac{1}{\ell}}\right)m(C).
\end{align*}
Choose  $\ell$ sufficiently large such that
\begin{equation}\label{mirival}
\frac{1}{2\varepsilon}\sum_{k\ge 0} \frac{(k+1)}{\ell}m(P(k))< \left(\frac{1}{2\varepsilon\log(R)}+ \frac{\beta_2}{4\varepsilon\log(R)}\right)V.
\end{equation}

Fix $n \in \mathbb{N}$ and $x \in \mathbb{S}^1$. Given $k \ge -\ell$, let us define 
\begin{equation}\label{100h}
\widetilde{P}_{n,k}(\omega,x):= \#\{i \in \{0,\dots, n-1 \} ; f^i_{\omega}(x) +\omega_i \in P(k)  \},
\end{equation}
analogously, consider
\begin{align*}
\widetilde{G}_n(\omega,x) &:= \#\{i \in \{0,\dots, n-1 \}; f^i_{\omega}(x)+\omega_i \in G  \}.
\end{align*}
Recall that $G= \{x\in\mathbb S^1; \|f'(x)\|> R\}$ and $\bigcup_{k \ge -\ell }P(k) = \mathbb{S}^1 \setminus G$. Therefore,
\begin{equation}\label{lowerb}
S_n(\omega,x):=\sum_{i=0}^{n-1}\log|f'_{\theta^i \omega}(f^i_{\omega}(x))| \ge \tilde{G}_{n}(\omega,x)\log(R) - \sum_{k =0}^{\infty}\widetilde{P}_{n,k}(\omega,x) \frac{(k+1)}{\ell}\log(R).
\end{equation}
Given $x \in \mathbb{S}^1$, consider the set of noise realisations 
\begin{align*}
A_n(x) := \left\{\omega\in\Omega; \  \widetilde{G}_n(\omega,x) - \sum_{k \ge 0} \widetilde{P}_{n,k}(\omega,x) \frac{k+1}{\ell} < \left(1-h-\left(\frac{1+\beta_2}{2\varepsilon \log(R)}\right)V\right)n \right\}.
\end{align*}
Note that by \eqref{zetauno} and  \eqref{lowerb} we have 
\begin{equation}\label{inclusion}
\left\{ S_n(\omega,x)< Z_1(h)n \right\} \subset A_n(x).
\end{equation}
To estimate the measure of $A_n(x)$, we introduce the set $A^1_n(x)$ of all noise realisations for which the orbit of $x$  enters one of the sets  $P(k)$ for  $k \ge n$ at least once,
\begin{align*}
A^1_n(x) &:= \left\{\omega\in \Omega;\ \text{there exists}\ k \ge n\ \text{such that}\  \widetilde{P}_{n,k}(\omega,x) >0  \right\}.
\end{align*}
 We further consider the event that the orbit of $x$ visits $G$ with a frequency smaller than  $1-h- \beta_2V/(8\varepsilon\log(R))$:
\begin{align*}
A^2_n(x) := \left\{\omega\in \Omega;\ \tilde{G}_{n}(\omega,x)\le \left(1-h-\frac{\beta_2 }{8\varepsilon \log(R)}V\right)n \right\}.
\end{align*}
The set $A_n(x)$ can be decomposed as 
\begin{equation}\label{decomposition}
A_n(x) := A^3_n(x) \sqcup \left( A_n(x) \cap \left(A^1_n(x) \cup A^2_n(x)\right)\right),
\end{equation}
where
\begin{align}
A^3_n(x) &=  A_n(x) \setminus \left(A^1_n(x)\cup A^2_n(x)\right)\nonumber\\
&= \left\{\omega \in \Omega;\  \sum_{k = 0}^{n-1} \widetilde{P}_{n,k}(\omega,x) \left(\frac{k+1}{\ell}\right)  > \widetilde{Z}n \right\}, \label{name}
\end{align}
and
\begin{equation}\label{name2}
\widetilde{Z} := \left(\frac{3\alpha}{8\varepsilon \log(R)}+\frac{1}{2\varepsilon \log(R)}\right)V.
\end{equation}
Observe that 
\begin{align}
\mathbb{P}(A_n(x)) \le \mathbb{P}\left(A^1_n(x)\right) + \mathbb{P}(A^2_n(x)) + \mathbb{P}(A^3_n(x)), \label{A123}
\end{align}
therefore to prove Lemma \ref{ledev}  it is enough to estimates on  the measure of $A^1_n(x)$, $A^2_n(x)$ and $A^3_n(x)$ this is done in Steps 2-4 below.
\begin{step}{2}\label{motivation}
We show that there exists a $D_0 \ge 1$, such that, for all $n \ge 0$,
 $\mathbb{P}\left(A^1_n(x)\right) \le D_0 R^{-\frac{n}{2\ell}}.$
\end{step} Let 
$\widetilde{P}(n):= \bigcup_{k \ge n}P(k)$, then
\begin{align*}
\mathbb{P}\left(A^1_n(x)\right) \le \sum_{i=0}^{n-1} \mathbb{P}\left( f^i_{\omega}(x)+x_i \in \widetilde{P}(n)  \right) \le \sum_{i=0}^{n-1} \sup_{y \in \mathbb{S}^1}\mathbb{P}\left( y+\omega_i \in \widetilde{P}(n)\right).
\end{align*}
Since the critical set of $f$ is non-degenerated, there exists a constant $D \ge 1$ such that ${m\left(\widetilde P(n)\right)\le D R^{-\frac{n}{\ell}}}$, for all $n \ge 0$. As a consequence, we obtain that
\begin{align*}
\mathbb{P}\left( y+\omega_i \in \widetilde{P}(n)  \right) \le \frac{1}{2\varepsilon}m(\widetilde{P}(n))
&\le D R^{-\frac{n}{\ell}},
\end{align*}
which yields $\mathbb{P}\left(A^1_n(x)\right) \le  DnR^{-\frac{n}{\ell}} \le D_0R^{-\frac{n}{2\ell}},$ for some $D_0 \ge 1$.
 This concludes Step 2. 
\begin{step}{3}\label{coollemma}
 There exists $q_1 \in (0,1)$ such that, for all $n \ge 0$, $\mathbb{P}\left(A^2_n(x)\right) \le q_1^n.$
\end{step}

Given $j\in \mathbb N,$ let us define
\begin{align}
\widehat{\pi}_j(A) = \mathbb{P}(\omega\in \Omega; (\omega_0,\dots,\omega_{j-1}) \in A) \qquad \ \text{for every} A \subset \mathcal{B}([-\varepsilon,\varepsilon]^j). \label{pihat}
\end{align}

Observe that $A^2_n(x)$ can be rewritten as
\begin{align*}
A^2_{n}(x)= \left\{ \sum_{i=0}^{n-1}\mathbbm 1_{\mathbb{S}^1 \setminus G}(f^i_{\omega}(x)+\omega_i)>\left(h+ \frac{\beta_2}{4}V\right)n \right\}.
\end{align*}
 By  the Markov inequality, for any $\zeta>1$, we obtain that
\begin{align*}
 \mathbb{P}\left(A^2_{n}(x)\right) \le \frac{\mathbb{E}\left[\zeta^{\sum_{i=0}^{n-1}\mathbbm 1_{\mathbb{S}^1\setminus G}\left(f^i_{\omega}(x)+\omega_i\right)}\right]}{\zeta^{\left(h+ \frac{\beta_2}{4}V\right)n}}.
\end{align*}
Therefore,
\begin{align*}
&\mathbb{E}\left[\zeta^{\sum\limits_{i=0}^{n-1}\mathbbm 1_{\mathbb{S}^1\setminus G}(f^i_{\omega}(x)+\omega_i)}\right]  = \int_{\Omega} \prod_{i=0}^{n-1} \zeta^{\mathbbm 1_{\mathbb{S}^1\setminus G}(f^i_{\omega}(x)+\omega_i)} \widehat{\pi} (\d \omega_0,\dots, \d \omega_{n-1}) \\
=& \int_{\Omega}\prod_{i=0}^{n-2} \zeta^{\mathbbm 1_{\mathbb{S}^1\setminus G}(f^i_{\omega}(x)+\omega_i)}\left[ \int_{[-\varepsilon,\varepsilon]} \zeta^{\mathbbm 1_{\mathbb{S}^1\setminus G}(f^{n-1}_{\omega}(x)+\omega_{n-1})}\mathbb P(\d \omega_{n-1})\right] \widehat{\pi}(\d \omega_0,\dots,\d \omega_{n-2}), 
\end{align*}

Note that, given $\omega_0,\dots,\omega_{n-2}$ fixed, we have that
$$
\int_{[-\varepsilon,\varepsilon]} \zeta^{\mathbbm 1_{\mathbb{S}^1\setminus G}(f^{n-1}_{\omega}(x)+\omega_{n-1})}\mathbb{P}(\d \omega_{n-1}) = (\zeta-1)\frac{\mathrm{Leb}(\omega_n ; f^{n-1}_{\omega}(x)+\omega_{n-1} \notin G)}{2\varepsilon}+ 1.
$$
Since  $\varepsilon > (1-m(G)/(2h)$, we obtain that  
\begin{align*}
\int_{[-\varepsilon,\varepsilon]} \zeta^{\mathbbm 1_{\mathbb{S}^1\setminus G}(f^{n-1}_{\omega}(x)+\omega_{n-1})}\mathbb P(\d \omega_{n-1}) \le \zeta h + 1-h.
\end{align*}
Repeating the argument recursively, we have that $\mathbb{E}\left[\zeta^{\sum_{i=0}^{n-1}1_{\mathbb{S}^1\setminus G}(f^i_{\omega}(x)+\omega_i)}\right]  \le (\zeta h + 1-h)^n.$ Hence, for any $\zeta>1$ 
\begin{align*}
 \mathbb{P}\left(A^2_n(x)\right)\le \left(\frac{h\zeta+ 1-h}{\zeta^{h+\frac{\beta_2}{4}V}}\right)^n.
\end{align*}
Taking $q_1:= (\zeta h+ 1-h)/(\zeta^{h+\frac{\beta_2}{4}V})$  for $\zeta$ close enough to $1$ we conclude the proof of Step 3. 
 
\begin{step}{4}\label{importandev}
We show that there exists  $q_2 \in (0,1)$ such that, for all $n \ge 0$, 
 $\mathbb{P}\left(A^3_n(x)\right) \le q_2^n,$ and conclude the proof of Lemma \ref{ledev}. 
\end{step}

By the Markov inequality and \eqref{name}, for all $1<\zeta$.
\begin{equation}\label{Chebishev}
\mathbb{P}(A^3_n(x)) \le \frac{\mathbb{E}\left[\zeta^{ \sum_{k = 0}^{n-1} \widetilde{P}_{n,k}(\omega,x)\frac{k+1}{\ell}}\right]}{\zeta^{n\widetilde{Z}}}
\end{equation}
Given $k \in \{0,\dots, n-1\}$ define $M_j(\omega,x) := \prod_{k=0}^{n-1}\zeta^{\frac{k+1}{\ell} 1_{P(k)}(f^j_{\omega}(x)+\omega_j)}$ and recall the definition of $\widehat{\pi}_k$ in \eqref{pihat}.

From \eqref{100h}
\begin{align*}
\mathbb{E}\left[\zeta^{ \sum\limits_{k = 0}^{n-1} \widetilde{P}_{n,k}(\omega,x)\frac{k+1}{\ell}}\right] &= \mathbb{E}\left[\prod_{k=0}^{n-1}\zeta^{\widetilde{P}_{n,k}(\omega,x)\frac{k+1}{\ell}}   \right]
= \mathbb{E}\left[\prod_{k=0}^{n-1}\prod_{j=0}^{n-1}\zeta^{\frac{k+1}{\ell}\mathbbm 1_{P(k)}(f^j_{\omega}(x)+\omega_j)} \right] \\
&=\mathbb{E}\left[\prod_{j=0}^{n-1}\left[\prod_{k=0}^{n-1} \zeta^{\frac{k+1}{\ell}\mathbbm 1_{P(k)}(f^j_{\omega}(x)+\omega_j)} \right]\right] 
= \int_{\Omega}  \prod_{j=0}^{n-1} M_j(\omega,x)\widehat{\pi}_n(\d \omega).
\end{align*}
By Fubini's theorem
\begin{equation}\label{fubini}
\int_{\Omega} \prod_{j=0}^{n-1}M_j(\omega,x)\widehat{\pi}_{n-1}(\d \omega)= \int_{\Omega}\prod_{j=0}^{n-2}{M_j(\omega,x)}\left[\int_{[-\varepsilon,\varepsilon]}M_{n-1}(\omega,x)\mathbb P(\d \omega_n)\right]\widehat{\pi}_{n-1}(\d \omega).
\end{equation}
Observe that, given $(\omega_0,\dots,\omega_{n-2}) \in [-\varepsilon,\varepsilon]^{n-1}$, we can rewrite
\begin{align*}
\int_{[-\varepsilon,\varepsilon]}M_{n-1}(\omega,x)\mathbb{P}(\d \omega_{n-1}) &= \int_{[-\varepsilon,\varepsilon]} \prod_{k=0}^{n-1} \zeta^{\frac{k+1}{\ell}\mathbbm 1_{P(k)}(f^{n-1}_{\omega}(x)+\omega_{n-1})}\mathbb{P}(\d \omega_{n-1})\\
&\le \max_{y \in \mathbb{S}^1}\int_{[-\varepsilon,\varepsilon]} \prod_{k=0}^{n-1} \zeta^{\frac{k+1}{\ell}\mathbbm 1_{P(k)}(y+\omega_{n-1})}\mathbb{P}(\d \omega_{n-1})
\\
&\le 1- \frac{1}{2\varepsilon} m\left(\bigcup_{k=0}^{n-1 }P(k)\right)+ \sum_{k=0}^{n-1} \zeta^{\frac{(k+1)}{\ell}}\frac{m\left(P(k)\right)}{2\varepsilon}.
\end{align*} 
As a consequence 
\begin{equation}\label{nickchata}
\int_{\Omega} \prod_{j=0}^{n-1}M_j(\omega,x)\widehat{\pi}_{n-1}(\d \omega)\le  W(\zeta) \int_{\Omega} \prod_{j=0}^{n-1}M_j(\omega,x)\pi_{n-2}(\d\omega),
\end{equation}
where  $W(\zeta):= 1+ \sum_{k=0}^{n-1} (\zeta^{\frac{(k+1)}{\ell}}-1)m\left(P(k)\right)/(2\varepsilon).$
Repeating the the same procedure $n-1$ times, we obtain
\begin{equation}\label{exponential}
\mathbb{P}(A^3_n(x)) \le \left(\frac{W(\zeta)}{\zeta^{\widetilde{Z}}}\right)^n. 
\end{equation}
Since $W(1)=1$ and, by \eqref{mirival} and \eqref{name2}
\begin{align*}
W'(1)&= \sum_{k=0}^{n-1}
\frac{(k+1)}{\ell}\frac{m\left(P(k)\right)}{2\varepsilon} < \left(\frac{1}{2\varepsilon \log(R)}+\frac{\beta_2}{4\varepsilon \log(R)}\right)V< \widetilde{Z}.
\end{align*}
Therefore $W(\zeta) < \zeta^{\widetilde{Z}}$, for $\zeta>1$ sufficiently close to $\zeta$. Therefore, $q_2:= W(\zeta)/\zeta^{\widetilde{Z}}<1,$ for $\zeta$ sufficiently close to $1$.

To conclude the proof of the lemma. Combining \eqref{exponential} with Steps 2 and 3, which show that \eqref{A123} decays exponentially fast which concludes the proof of the lemma.
\end{proof}
 
Finally, we prove Proposition \ref{LHiperbolicTijmens}.

\begin{proof}[Proof of Proposition \ref{LHiperbolicTijmens}] {Let $c>0$ be fixed, so that $\varepsilon \ge \alpha^{c-1}$.  Choose $\alpha_1> 1$ sufficiently large such that for every $\alpha > \alpha_1,$  the set $G(\alpha):=\{x\in\mathbb S^1;\ \|d f(x)\|> \alpha^{c/2}\},$ satisfies the existence of an interval $I\subset G$ such that $f(I) = \mathbb S^1$. Moreover, define $s = s(\alpha):= m(G(\alpha)).$ Consider $C = \{x\in\mathbb S^1;\ |f'(x)|\leq 1\}.$ From Lemma \ref{LGiu} and \eqref{LGiu2}, we can choose $\alpha_2 > \alpha_1$, such that for every $\alpha>\alpha_1$
$$  h = \frac{1}{2}  < 1 + \frac{2\alpha^{1-c }}{\log(\alpha^{c/2})} \int_C \log |f'(x)| \d x\ \text{and}\  \frac{1-s(\alpha)}{h} < \frac{1}{\alpha^{1-c}}.$$

For every $\sigma^2, r>0$, $x\in\mathbb S^1$ and $n\in \mathbb N$ define
\begin{align*}
A(n,x,r)&:= \left\{\omega \in \Omega;\sum_{i=0}^{n-1}\log(\mathrm{dist}_{r}(f^i_{\omega}(x),\mathscr{C}_{\theta^j \omega}))<r n\right\},\\
B(n,x,\sigma^2)&:= \left\{\omega \in \Omega;\sum_{i=0}^{n-1}\log| f'_{\theta^i \omega}(f^i_{\omega}(x))|\le \sigma^2 n \right\}.
\end{align*}

From Lemmas \ref{largedeviazioni} and  \ref{ledev}, we obtain that there exists $r$ and $\sigma^2$ small enough such that, for every $\alpha > \alpha_1,$
\begin{align}
\sup_{x \in \mathbb{S}^1}\mathbb P [ A(n,x,r) \cup B(n,x,\sigma^2)]\xrightarrow[]{n \to \infty} 0,\label{Giugiuexp}
\end{align}

From the Pliss lemma (see \cite[Lemma 6.2]{Alves}) and the definition of $L$-sparse hyperbolic times, provided  $r$ and $\sigma^2$ are small enough, there exists $\gamma \in (0,1)$ such that for every $x\in\mathbb S^1,$
\begin{equation}\label{eq1}
\left\{\omega\in\Omega;\#\left(\{1,\ldots,n\} \cap \{\tau_i(\omega,x)\}_{i=1}^{\infty}\right) \leq   \frac{\gamma}{2 (L+1)}  n\right\} \subset  A(n,x,r) \cup B(n,x,\sigma^2).
\end{equation}
For more details on the construction of $\gamma$, see \cite[Proposition 6.3]{Alves}. The proposition follows by combining equations \eqref{Giugiuexp} and \eqref{eq1}.}

\end{proof}

\bibliographystyle{plain} 
\bibliography{main}
\end{document}